\documentclass[11pt]{article}

\title{Invariance principle for `push' tagged particles for a Toom Interface}

\author{Nicholas Crawford and Wojciech De Roeck}
\usepackage{amsmath, amsfonts, amssymb, amsthm}
\usepackage{enumerate, color}
\usepackage{multido}
\usepackage{url}
\usepackage{mathrsfs}
\usepackage{mathtools}
\usepackage{verbatim}
\usepackage{setspace}
 \usepackage{hyperref}
\usepackage{cleveref}
\usepackage{cite}

\usepackage[margin=3cm]{geometry}
\usepackage{tikz}

\newtheorem{theorem}{Theorem}[section]
\newtheorem{lemma}[theorem]{Lemma}
\newtheorem{remark}[theorem]{Remark}
\newtheorem{proposition}[theorem]{Proposition}
\newtheorem{definition}[theorem]{Definition}
\newtheorem{assumption}[theorem]{Assumption}
\newtheorem{corollary}[theorem]{Corollary}

\crefname{theorem}{Theorem}{Theorems}
\crefname{lemma}{Lemma}{Lemmas}
\crefname{proposition}{Proposition}{Propositions}
\crefname{section}{Section}{Sections}
\crefname{assumption}{Assumption}{Assumptions}
\crefname{equation}{Eq.}{Eqs.}

\newcommand{\R}{\mathbf{R}}
\newcommand{\N}{\mathbf{N}}

\newcommand{\E}{\mathbf{E}}
\newcommand{\Z}{\mathbf{Z}}

\renewcommand{\Pr}{\mathbf{P}}

\newcommand{\Om}{\Omega}
\newcommand{\Bp}{\textrm{Ber}_p}

\renewcommand{\d}[1]{\nabla_{#1}}

\newcommand{\eps}{\epsilon}

\newcommand\otimesal{\mathop{\hbox{\raise 1.6 ex
  \hbox{$\scriptscriptstyle\mathrm{al}$}
\kern -0.92 em \hbox{$\otimes$}}}}
\newcommand\oplusal{\mathop{\hbox{\raise 1.6 ex
  \hbox{$\scriptscriptstyle\mathrm{al}$}
\kern -0.92 em \hbox{$\oplus$}}}}
\newcommand\Gammal{\hbox{\raise 1.7 ex
\hbox{$\scriptscriptstyle\mathrm{al}$}\kern -0.50 em $\Gamma$}}

\let\ka=\kappa   
\let\si=\sigma

   \let\Om=\Omega
  \let\Si=\Sigma

\newcommand{\caA}{{\mathcal A}}
\newcommand{\caB}{{\mathcal B}}
\newcommand{\caC}{{\mathcal C}}

\newcommand{\caF}{{\mathcal F}}
\newcommand{\caG}{{\mathcal G}}

\newcommand{\caJ}{{\mathcal J}}

\newcommand{\caM}{{\mathcal M}}

\newcommand{\caO}{{\mathcal O}}

\newcommand{\caS}{{\mathcal S}}
\newcommand{\caT}{{\mathcal T}}

\newcommand{\caW}{{\mathcal W}}

\newcommand{\bbN}{{\mathbb N}}

\newcommand{\opunit}{\text{1}\kern-0.22em\text{l}}

\newcommand{\cf}{cf.\;}

\newcommand{\norm}{ \|}
\newcommand{\str}{ |}

\newcommand{\e}{{\mathrm e}}

\renewcommand{\d}{{\mathrm d}}

\newcommand{\beq}{ \begin{equation} }
\newcommand{\beqs}{ \begin{equation*} }
\newcommand{\eeqs}{ \end{equation*} }
\newcommand{\eeq}{ \end{equation} }
\newcommand{\bet}{ \begin{theorem} }
\newcommand{\eet}{ \end{theorem} }

\newcommand{\prob}{\mathbf P}

\begin{document}
\maketitle
\begin{abstract}
In many interacting particle systems, tagged particles move diffusively upon subtracting a drift.  General techniques to prove such `invariance principles' are available for reversible processes (Kipnis-Varadhan) and for non-reversible processes in dimension $d>2$. The interest of our paper is that it considers a non-reversible one-dimensional process: the Toom model. The reason that we can prove the invariance principle is that in this model, push-tagged particles move manifestly slower than second-class particles.

\end{abstract}

\section{Introduction}\label{S:Intro}

Let us introduce the Toom model. 
It plays on spin configurations $\sigma:=(\sigma({x}))_{x \in \Z} \in \Omega$ with $\Omega= \{-1, 1\}^\Z$, but it is good to think of the different values of $\sigma(x)$ as the site $x$ being occupied by either $+$ or $-$ particles.   Each $\pm 1$ particle is equipped with an exponential rate $\lambda_{\pm}$ clock.  When the clock rings for a particle of sign $\eta$, the particle exchanges positions with the first particle to its right of opposite sign $-\eta$. Since that opposite sign particle can be arbitrarily far away, this process is of infinite range, it is not a Feller process.  Here and below, we'll refer to this process as $\sigma_t:= (\sigma_{t}(x))_{x \in \Z}$.  The Bernoulli measures $\Bp$, where $p={\Bp}(\sigma(x)=+)$, are invariant, and, in fact, we have showed \cite{crawford2015toom} that they are the only invariant measures satisfying certain regularity conditions.  In what follows, we are always referring to these stationary processes.

%One remarkable feature of this model is that the restriction from $\Om[\N]$ to the first $N$ vertices $\Om[N]:=\{-1, 1\}^{[N]}$ is itself a Markov chain; the dynamics is the same unless a clock rings for a spin in the last block of constant sign in $[N]$.  For updates of spins in the last block, the dynamics reduces to single vertex spin flips.  As a result, the restricted chain is irreducible on $\Om[N]$ and has a unique stationary measure $\pi_N$.  As the sequence of measures $(\pi_{N})_{N\in \N}$ is consistent, this in term implies that the full chain has a unique invariant measure on $\Om[\N]$, $\pi_{\infty}$, which restricts to $\pi_N$ on $\Om[N]$. 

In the above description, there is an obvious notion of 'tagged particle', but it is not this notion for which we can prove the invariance principle. Instead, we consider Push-tagged particles:    Let's focus on a single signed particle and suppose that the block of spins to its immediate right has the same sign as the particle. Then, rather than viewing the particle as jumping \textit{over} its neighboring block to the right,  we can view the particle as moving to its right one site.  In doing so it \textit{pushes} the entire right neighboring block of particles one site right as well. 

 It is clear that this dynamics leads to the same \textit{unlabeled} particle system as was defined above.   This description of the dynamics provided inspiration for the paper \cite{borodin2007large} in which the authors discussed a model called  Push-ASEP which has a integrable structure.   Note that Push-ASEP is really just the "totally asymmetric" case of the present setup with $\lambda_+=1, \lambda_-=0$, see also a generalization: q-pushASEP,  in \cite{corwin2015q}

A convenient feature of this "pushing" description is that the dynamics preserve any total ordering of particles of the same type.  That is, if we denote by $Y^{(j)}_t$ the position, at time $t$, of the particle which starts at $j$ and if $x<y$, then on the event that $\sigma_{0}(x)=\sigma_0(y)$
\[
Y^{(x)}(t) \leq Y^{(y)}(t) \text{ for all $t>0$.}
\]
This simple observation is important in our proofs. 
Our main result is 
\begin{theorem}[Functional CLT for Tagged Push Particles]
\label{T:CLTtag}
Fix $\lambda_{+},\lambda_- $ nonnegative not both zero, and $p \in (0, 1)$. 
Starting from the Bernoulli measure $\Bp$ conditioned on $\sigma_0(0)=\pm 1$ (i.e.\ fixing the sign of the push particle), 
\[
 \frac{Y^{(0)}_{ nt}- v_{\pm} n t }{\sqrt n}\quad \mathop{\Rightarrow}^{d} \quad \sqrt{D} B_t,\qquad  t\in [0,1]
\] 
where $B_t$ is standard Brownian motion, the convergence is in distribution on Skorohod space and the drift is given by
\[
v_{\pm}=\lambda_{\pm}\left(\frac{1}{1-p}\right)-\lambda_{\mp}\left(\frac{1-p}{p^2}\right)
\]
The diffusion constant is positive $D>0$. 
\end{theorem}
%
%Just as for conventional particles, the environment seen from the push-particle is also $\Bp$ and stationary. 
%

%allows us to apply ideas from the proof of \cref{T:CLT} to prove a functional CLT for $Y^{(0)}_t$.
%
% It is easy to see that such defined 'push-particles' of the same sign preserve their ordering. 

%
%Let $Y_t$ be the position of the push-particle at time $t$, and let us for simplicity assume that it is a $+$ particle. When the clock of this particle rings, we push it forward one-site $Y_t=Y_{t-}+1$, regardless of how far the conventionally defined $+$ particle has jumped forwards.  When the clock rings on a $+$ particle to the left of $Y_{t-}$, and this clock rings brings that particle in one jump to the right of $Y_{t-}$, then we also put $Y_t=Y_{t-}+1$. When the clock rings on a $-$ particle at $x$ to the left of $Y_{t-}$, such that  this clock rings brings onto $Y_{t-}$, then we put $Y_{t}=x$, i.e.\ the particles exchange position.  
%
%\wdr{Or the old version?::}
%
%
%
%\wdr{end of insertion}

%Let us call $Y_t$ the position of the push-particle initially started from time $0$
%
%Our main result is hence 
%$$
%\frac{1}{\sqrt{n}}(Y_{nt}-n\E(Y_t))  \Rightarrow B_t, \qquad t \in [0,1] 
%$$
The technique yields at the same time invariance principles for additive functions like 
$$
X_t= \int_{0}^t \textrm{d} s \: \sigma_s(0).
$$
and integrated currents, like the total number of $\pm$ particles crossing a given edge. However, for the additive functionals, we do not prove positivity of the variance.

Let us conclude by reviewing some earlier results on functional CLTs for tagged particles in conservative particle systems.  The classic paper by 
Kipnis-Varadhen,
\cite{kipnis1986central}, implies CLTs for symmetric exclusion processes (except the nearest neighber case) while \cite{Varadhan} extends this to general zero drift jump kernels.  Both results work in any dimension. For non-zero drift, there is a general approach \cite{sethuraman2000diffusive}
for asymmetric exclusion processes in dimension $d\geq 3$.  In dimensions $d=1, 2$, there is no general approach available and results can only be proved on a case by case basis using specific features of the underlying models.  Moreover, results in this case seem to be few and far between see for example
\cite{RWRW,Kipnis}.  {Our result also uses specific properties of the model, in particular fast mixing exhibited via natural coupling and the order-preserving feature of push particles.}

\subsection{Preliminaries: the Dynamics on $\Z$.}
We fix once and for all $\lambda_{\pm}>0$ with $\lambda_++\lambda_-=1$, which just sets the overall time scale. 
We consider a sequence of i.i.d.\ rate one Poisson point processes $(N_x(t))_{x\in \Z}$ associated with vertices $x\in \Z$.  Besides these Poisson point processes, the sample space on which our processes are defined supports a two dimensional array of of i.i.d.\ uniform $[0,1]$ variables $(U_{x, j})_{x\in \Z, j\in \N}$.  Let $(\varOmega, \mathbb P; \caB_{\varOmega})$ denote a probability space which supports all these variables. Define the filtration of sigma algebras $(\caF_t)_{t \in \R^+}$ on $\caB_{\varOmega}$ by
\[
\caF_t = \sigma\left(N_x(s): s \leq t; \: U_{x, k}: k \leq N_x(t)\right).
\]
Let us define $\Om=\{-1,+1\}^{\Z}$ and equip $\Om$ with with its natural product topology and associated Borel sigma algebra $\caB$.  Let $\Sigma= \Om\times \varOmega$ and equip it with its natural product sigma algebra.  
Finally, let $I\subset \R_+$ be any closed interval and let
 $D_{\Om}(I)= D(I\to\Om)$ be the space of c\'adl\'ag functions from $I$ to $\Om$. In case $I=[0,\infty)$ we simply denote this space by $D$.  We equip $D_{\Om}(I)$ with the Skorokhod topology and associated Borel sigma algebra, the latter being denoted $\caB(D_\Omega([0,\tau]))$.

In general, given a pair of measurable spaces $(\mathbf X, \caF); (\mathbf Y, \caG)$ and a family of random variables $(X_i)_{i \in I}$, we shall denote by $\caB(X_i: i\in I)$ the sigma algebra generated by the $X_i's$.  Also, given a measure $\mu$ on $(\mathbf X, \caF)$, the Lebesgue space $L^q(\mathbf X, \mu), q\geq 1$ will often be abbreviated $L^q(\mu)$, and even $L^q$ when confusion is unlikely, with corresponding norm denoted by $\|\cdot\|_{L^q(\mu)}$.

 As already remarked, the process is non-Feller and therefore cannot be defined in the standard way (see \cite{crawford2015toom} for further discussion on this point). Nevertheless, in   \cite{crawford2015toom} we constructed the process  with $\Bp$ initial conditions.   For $\lambda_\pm$ fixed, it is convenient to introduce the thinned Poisson processes $N_{x, \pm}(t)$ by the differentials $\d N_{x,+}(t)= \mathbf{1}\{U_{x, N_x(t)}< \lambda_{+}\}\d N_{x}(t)$ and $\d N_{x,-}(t)= \d N_{x}(t)-\d N_{x,+}(t)$.
 \begin{theorem}\label{def: existence dynamics}
  There is a $\mathbf{P}_{\Bp}$-a.s. defined 
random variable $F: \Sigma \rightarrow D$, i.e.\ a c\'adl\'ag process,  such that if we denote the value of $F$ at time $t$ by $\si_t$
\begin{enumerate}
\item (Stationarity) $\si_t$ is $\Bp$-distributed for any $t$.
 \item (SDE is satisfied)
{The SDEs
\begin{multline}   \label{eq: infinite l sde}
 \si_{t_2}(x)-\si_{t_1}(x) =\sum_{\eta =\pm 1} -2\eta \int_{t_1}^{t_2} \d t \chi^{\eta}_{x}(\si_{t_-}) \,  \d N_{x,\eta}(t)\\
 + \sum_{\eta =\pm 1} 2\eta \int_{t_1}^{t_2} \d t \sum_{ y < x}  \chi^{\eta}_{[y,x-1]}(\si_{t_-}) \chi^{-\eta}_{x}(\si_{t_-}) \,  \d N_{y,\eta}(t)
\end{multline}
are satisfied $\mathbf{P}_{\Bp}$-a.s.  In particular, the right hand side is absolutely summable and that the equality \eqref{eq: infinite l sde} holds for any $x$ and $t_1<t_2$.}
%
% \item  (Finite speed of information propagation).  For any $T>0$ and any finite set $S \subset \Z$,
% \beq
% \label{E:FV}
% \sup_{t \in [0,T), x \in S} \str \si_t(x)- \si^L_t(x)\str \qquad \text{converges to $0$ in distribution as $L\to \infty$.}
% \eeq 
 \end{enumerate}
\end{theorem}
\begin{proof}
This is a restatement, in slightly different language, of Lemma 2.8 of \cite{crawford2015toom}.   The main point here is the language of stochastic differentials to describe the evolution of spin variables.
\end{proof}

\begin{remark}
\label{R:Adj}
The $(U_{x, k})_{x, k}$ here may seem obscure.   Through these variables we can couple an arbitrary collection of Toom trajectories $(\si^{j}_t)_{j\in J}$ indexed by an at-most countably infinite index set $J$.  
This was used in a variety of ways in \cite{crawford2015toom}.  In particular, we recall their use in \Cref{def: existence dynamics}. 
 The key, and most concrete, step of the proof of \Cref{def: existence dynamics} was the fact that, for short times one can couple a sequence of finite systems $\si_t^L$  with periodic boundary conditions (on $[-L, L)$ say) so that for each finite window $[-K, K)$ and all $t\in[0, \eps)$, $\lim_{L}\si^L_t$ exists $\mathbf{P}$-a.s. and is $\Bp$ distributed for all $t\in[0, \eps)$.  We will need this fact in \cref{S:Adjoint} to verify a time reversal identity between a Toom process moving to the right and a Toom process moving to the left.
\end{remark}

In this paper, we only consider couplings between $\si^{j}_t$'s whose initial distribution is $\Bp$.
Formally,
\begin{definition}\label{def:coupling}
Let $\{(\si^j_t):j\in J\}$ be two or more Toom processes (not necessarily on the same subset of $\Z$) having respective initial distributions $\Bp$. When we discuss a ``coupling of the $\{(\si^j_t):j\in J\}$ started from $\mu$'' we mean the following: $\mu$ is assumed to be a measure on $\prod_{j\in J}\{\pm 1\}^\Z$ whose marginals are the $\Bp$. The coupling is then the collection of the $D$-valued random variables $\si^j$ given by
\[
\sigma^j=F(\eta^j,\omega)\qquad\sigma^j:\Big(\prod_{j\in I}\{\pm 1\}^\Z\Big)\times\varOmega\to D.
\]
\end{definition}
The existence of a coupling started from a given $\mu$ is immediate from the fact that $F$ is a $\mathbf{P}_{\Bp}$-a.s. almost-surely defined function and each single-spin-configuration marginal of $\mu$ is $\Bp$.  The law of $\{(\si^j_t):j\in J\}$ starting from an initial measure $\mu$ will be denoted by $\mathbf{P}_{\mu}$.

\subsection*{Acknowledgements}
%We are grateful to Gady Kozma for hel
%This work is a spin-off of the paper \cite{crawford2015toom} that was written at the same time, together with Gady Kozma. We are 
  WDR acknowledges the support of the DFG (German Research Fund), the Belgian Interuniversity Attraction Pole  P07/18 (Dygest) and the FWO (Flemish Research Council).  NC is supported by Israel Science Foundation grant number 915/12.

\section{Statement of Results and the Key Lemma}

Our method of proof is fairly flexible, applying in wider generality than indicated in \Cref{S:Intro} with little extra overhead.  It is convenient for us to formalize the collection of observables which satisfy functional CLTs. Let us consider processes $X(t)$ of the following form. 
\beq
\label{E:Semi}
\textrm{d} X(s)= \langle g(\si_{s-}) , \d N(s) \rangle + f(\si_{s-}) \textrm{d} s.
\eeq
where we used the notation
\[
 \langle g(\si_{s-}) , \d N(s) \rangle= \sum_{\eta, x}  g_{\eta, x} (\si_{s-}) \d N_{\eta, x}(s)
\]
and $f,g_{\eta,x}$  are measurable functions $ \Om  \rightarrow \R$. 

We need a few assumptions on $f,g_{\eta,x}$. The first one imposes some regularity, in particular implying that equation \Cref{E:Semi} is well-defined. 
\begin{assumption}[Finite $L^q$-norms] \label{ass: norms}
For any $1\leq q <\infty$, 
$$
\norm f \norm_{L^q(\Bp)} +  \sum_{\eta, x} \norm g_{x,\eta} \norm_{L^q(\Bp)}   <\infty
$$
\end{assumption}
This assumption implies that $X(t)$ has finite variation on any finite interval, \cf \Cref{def: existence dynamics} $\mathbf P$ a.s. 
Furthermore, we shall assume that the functions $f,g_{x,\eta}$ are well-approximated by local functions. For $f \in L^1(\Bp)$, we consider the conditional expectations
$$
 \texttt{P}_R f(\sigma) := \E_{\Bp}[f \, \big\str \,  \sigma(x), x \in [-R,R] ].
$$
 Our second assumption reads 
 \begin{assumption}[Local approximation] \label{ass: local}   
$$
\norm f-\texttt{P}_Rf \norm_{L^2(\Bp)} + \sum_x   \norm g_{x,\eta}-1_{\str x\str \leq R }\texttt{P}_Rg_{x,\eta}\norm_{L^2(\Bp)}   < C \e^{-c R}
$$
\end{assumption}
Note that  \Cref{ass: local,ass: norms} imply the bound in \Cref{ass: local} holds with $L^2$ replaced by $L^q, q\geq 1$ Finally, we give a condition which restricts $X(t)$ to be  measurable w.r.t.\ the path $\sigma_t$, i.e.\ to not depend on arrivals of the processes $N_{x,\eta}$ that have no bearing on the path $\sigma_t$. 
\begin{assumption}[Path measurability]  \label{ass: real}
For any $x,\eta$,
$$
g_{x,\eta}(\sigma) =  \chi(\sigma(x)=\eta) g_{x,\eta}(\sigma) 
$$
\end{assumption}

Let $Y_t$ denote the position of a tagged $+$-particle with initial position $0$.  The process $(\si_t, Y_t)$ with state space $\{(\si, y\in \Om \times \Z: \si(y)=1\}$ is Markovian.  We denote by $\mathbf P_{\Bp, y}$ the probability measure for this process where $\si_t$ is started from the measure $\Bp$ conditioned on the presence of a $+$-particle at $y$ and $Y_0=y$.   If we denote the spatial shifts $\tau_y:\Omega \rightarrow \Omega$ by $(\tau_y \si)_x:= \si_{y-x}$, this process descends via the map $(\si, y) \mapsto \tau_y \si$ to a Markov process on $\{-1, 1\}^{\Z \backslash \{0\}}$ called the \textit{environment seen from the (push) particle}.
It is easy to check the following, which is crucial for some of our results.
\begin{lemma}
\label{L:PushBer}
The Bernoulli measures on $\{-1, 1\}^{\Z \backslash \{0\}}$ are stationary when we pass to the environment-seen-from-the-push-particle perspective.
\end{lemma}
We will also consider CLTs for processes $X(t)$ defined by the equation
\beq
\label{E:Semitag}
\textrm{d} X(s)= \langle  \tau_{Y_{s-}} g(\si_{s-}) ,\tau_{Y_{s-}} \d N(s) \rangle +  \tau_{Y_{s-}} f( \si_{s-}) \textrm{d} s
\eeq
 $\tau_{y} \d N_{x, \eta}(s):= \d N_{y+x, \eta}(s)$ and 
\[
\tau_{y} \cdot g_{x, \eta}(\si):=g_{x, \eta}(\tau_{-y} \cdot \si).
\]
where $f, g$ satisfy the assumptions above.

In general, we will refer to processes $X(t)$ defined by  \Cref{E:Semi} or \Cref{E:Semitag} with $f,g$ satisfying the three Assumptions above as \emph{quasi-local processes}. If it is necessary to distinguish between \Cref{E:Semi} and \Cref{E:Semitag}, we will refer to the latter as 'quasi-local w.r.t. to the tagged particle' and the former as quasi-local w.r.t. to the origin.  We will often drop the subscript $\mathbf P_{\Bp}, \mathbf P_{\Bp, y}$ from our expressions below when there is no danger of confusion.  One exception to this is the exposition of \Cref{S:D}, where we deal with a coupling process and various initial measures.

We are now ready to state the main result. 
\begin{theorem}\label{thm: main}
Any quasi-local process $X(t)$ as defined above satisfies a Brownian invariance principle, i.e.\ the sequence of processes
$$
  \frac{1}{\sqrt{n}} (X({nt}) - tnv_X)), \qquad   0\leq t \leq 1,  \qquad  
$$
(with  drift $v_X :=(1/nt) \E(X({nt}))  <\infty$) converges weakly, as $n \to \infty$, to a multiple of Brownian motion, in the Skorohod topology. 
\end{theorem}
We list three important examples of such processes $X(t)$:
\begin{corollary}
In particular, the invariance principle holds for 
\begin{enumerate}
\item  Additive functionals  
$$
X(t)= \int_0^t \d s f(\sigma_{s}) 
$$
with $f$ satisfying the localization assumption.
\item  Tagged push-particles
$$
X(t)=Y_{t}.
$$
\item Integrated $\eta$-particle currents from $(-\infty, x)$ to $[x,\infty)$: 
\beq \label{eq: instant current 1}
X(t)=\sum_{y<x} \int_0^t   \chi^{\eta}_{[y,x-1]} (\si_{s-})  \d N_{y,\eta}(s).
\eeq
\end{enumerate}
\end{corollary}
To establish this corollary, we should check that these processes are indeed quasi-local process in the sense outlined above. For the first and third example, this is obvious, so we only comment on the   tagged particle, $Y_t$.  Let us consider only the case of the $+$-particle.  In the other case, there is an analogous representation.
Let 
\begin{align}
&Q_{x,r}(y, \sigma):=\mathbf 1\{y \in [x, x+r)\} \chi^+_{[x, x+r)}\chi^-_{\{x+r\}},\\
&P_{x,r}( y, \sigma) := \mathbf 1\{y = x+r\} \chi^-_{[x, x+r)}\chi^+_{\{x+r\}},
%&Q_{xr, s}= Q_{xr}(s, Y_{s-}, \sigma_{s-})\quad P_{xr, s}=P(s, Y_{s-}, \sigma_{s-}).
\end{align}
with empty products are treated as $1$. One can check that 
\[
{Y}_t=\sum_{x\in \Z}  \sum_{r>0} \int_0^t  Q_{x,r}(Y_{s-}, \sigma_{s-})\textrm{d} N_{+,x}(s) - \sum_{r>0}\int_0^t r P_{x,r}(Y_{s-}, \sigma_{s-}) \textrm{d} N_{-,x}(s)
\]
Note here that by \Cref{L:PushBer} the drift of $Y_t$ satisfies
\[
\E_{\Bp, 0}[{Y}_t]= v_{Y} t
\]
with
\[
v_Y:=  v_Y(\lambda+, \lambda_-, p)=\lambda_+\left(\frac{1}{1-p}\right)- \lambda_-\left(\frac{1-p}{p^2}\right).
\]
This establishes our main result Theorem \ref{T:CLTtag}, except for the positivity of the diffusion constant, which is however clear from the representation given in Section \ref{S:CLT}, where manifestly $D_1>0$ and $D_2 \geq 0$. 

\subsection{The Key Lemma}

The perspective we shall take in proving our results is that $\textrm{d} X(t)$ is a random signed measure on any finite interval $I \subset \R^+$.  Indeed, any real function of bounded variation defines a finite signed Borel measure. As remarked above, $X$ is  indeed a.s.\ of bounded variation on finite intervals.  The space of finite signed Borel measures over a compact set ${\cal X} \subset \R^d$, equipped with the total variation norm, is a Banach space that we denote by $\caM({\cal X} )$. It is the dual of $\caC_b({\cal X} )$, the bounded continuous functions on $X$ with the supremum norm. In all what follows, we take ${\cal X} $ some finite rectangle in $\R^d$.  Adding some standard considerations on Skorohod topology, we then derive
\begin{lemma}\label{lem: def measures}
Fix a finite interval $I$. On $(\Si,\mathbf P)$, we have almost surely defined random variables $\nu$, taking values in  $\caM(I)$, and given by $\nu(\d t):=\d X(t)$.
\end{lemma}

%\begin{proof}
%Only the measurability of the map yielding the random measures $\nu$ has not yet been explained. It suffices to show that $X(t) \mapsto \d X(t)$ as a map from Skorohod space $D_{\Om}(I)$ to  $\caM(I)$ is measurable. Moreover, since the pre-dual $\caC_b(I)$ is separable, the measurable structure on $\caM(I)$ (generated by its norm-topology) is the same as that generated by its weak$^*$-topology, i.e.\ the topology of weak convergence of measures. 
%Therefore it suffices to show measurability of $X(t) \mapsto \d X(t)$ where the paths are equipped with Skorohod topology and the measures are equipped with the topology of weak convergence. However, the mapping is continuous when $\caM(I)$ is given its weak* topology and hence is a measurable mapping. 
%\end{proof}

Whenever we consider expressions involving multiple quasi-local processes, we will index them as $X^{(i)}$, with corresponding integrands denoted by $f^{(i)}, g^{(i)}$.  To alleviate possible confusion, let us explicitly remark here that we will never mix the two cases of quasi-local processes and quasi-local w.r.t. a tagged particle. 
Now, given a finite collection $(X^{(i)})_{i=1}^{\ell}$ of  quasi-local processes,  the (random) product measure $\prod_i \textrm{d} X^{(1)}(s_i)$ is defined on the hypercube $[0,L]^{\ell}$ and we will always restrict these measures to the open simplex
$$
\caS_{\ell}(L) = \{ (t_1, \ldots, t_\ell) \in [0,L]^\ell,  t_j <t_{j+1}\},
$$

Most of our bounds will be phrased in terms of the variation of such measures. In particular, the key technical lemma we shall prove in the paper is stated as follows. 
 Let $T> 1, \kappa > 1$  and set, for $l =1,\ldots, \ell-1$
\beq
\label{E:EE}
E_l(T, \kappa)= \{(t_1, \ldots, t_\ell) \in \caS_{\ell}(L): \: t_l-t_1 \leq T,   t_{l+1}-t_l  \geq T^{\kappa}\}.
\eeq
Constants are allowed to depend on the processes $X^{(i)}$, in particular on the $f^{(i)}, g^{(i)}$, unless explicitly stated otherwise. 
\begin{lemma}
\label{L:Factor1}
Let $(X^{(i)})_{i=1}^{\ell}$ be  quasi-local processes. Let $\mu$ be the measure given by
\beq
\label{E:M1}
\mu(\d t_{1, \dotsc,\ell}):= \E\left[  \textrm{d} X^{(1)}(t_1) \cdots \textrm{d} X^{(\ell)}(t_\ell)\right] 
- \E\left[   \textrm{d} X^{(1)}(t_1) \cdots \textrm{d}X^{(l)}(t_{l})\right] \E\left[   \textrm{d} X^{(l+1)}(t_{l+1}) \cdots \textrm{d} X^{(\ell)}(t_\ell)\right].
\eeq
Then its variation $\str \mu \str$ on $E_l(T, \kappa)$ satisfies the bound
$$
|\mu|(E_l(T, \kappa)) \leq C L^{\ell-l+1} \e^{-c T^{\frac{\kappa-1}{2}}}.
$$
\end{lemma}
This lemma sets the stage for us to prove functional CLTs via the method of moments, see \Cref{S:CLT}.

\section{Bounds on iterated integrals and random measures}
\label{S:DofC}
This section provides a-priori bounds on the total variation of measures of the form
\beq \label{def: nu}
\nu(\d t_{1,\ldots, \ell}) =  \textrm{d} X^{(1)}(t_1) \cdots \textrm{d} X^{(\ell)}(t_\ell)
\eeq
on $\caS_{\ell}(L)$.  For a {quasi-local process} $X$, we use $X_R$ to denote the local approximation to $X$ obtained by replacing  $f$ by $f_R:= \texttt{P}_R f$ and $g_{x,\eta} $ by $g_{R, x, \eta}:=  1_{\str x \str \leq R}\texttt{P}_R g_{x,\eta}$ (see \Cref{ass: local}). 
% \red{\textbf{I dont understand (the function of) this sentence.  May just delete?}In what follows, we will not need to keep track of the norms of $f,g$ except in the cases where we are considering the above bounds for $f-\texttt{P}_Rf$ and $g-1_{\str x \str \leq R}\texttt{P}_Rg$.  }
%\wdr{This is about when yes or no one is keeping track of these functions in error bounds. Perhaps indeed unclear.}

{\begin{lemma}[A priori bounds]\label{lem: moment for use}
Let $\nu$ be the measure defined in \Cref{def: nu} and let $\nu_R$  be the same but with all $X^{(i)}$ replaced by $X^{(i)}_R$.   Then, for any $R$,
$$
\norm \str \nu \str ([0,t]^{\ell}) \norm_{L^q(\mathbf P)} \leq  C(1+t^{\ell}),  \qquad \norm \str \nu-\nu_R \str ([0,t]^{\ell}) \norm_{L^q(\mathbf P)} \leq  C(1+t^{\ell})\e^{-c R}
$$
with $C$ depending on $f,g$ but not on $t,R$. 
\end{lemma}}

To prove this lemma, the basic strategy will be to bound
$$
\str \nu  \str (\caS_{\ell}(L))   \leq \prod_{j=1}^\ell  |\nu^{(j)}|([0,L])
 $$
where  $\nu^{(j)}(\d t)=\textrm{d} X^{(j)}(t)$.  One can apply H\"older's inequality to the RHS to get
$$\E[
\norm \str \nu \str (\caS_{\ell}(L)) \norm_{L^q(\mathbf P)} \leq \prod_{j=1}^\ell   \left\norm\str\nu^{(j)}\str([0,L]) \right\norm_{L^{q\ell}(\mathbf P)}
$$
Note that 
$$
\str\nu^{(j)}\str([0,L]) = \bar X^{(j)}(L)- \bar X^{(j)}(0)
$$
where $\bar X^{j}(t)$ is obtained from $ X^j(t)$ by replacing $f,g$ with $\str f \str, \str g \str$ in the definition of the process $X$. Obviously, $\bar X^{(j)}(t)$ is a quasi-local process and hence our task reduces to proving bounds on $L^k$ norms of 
$$
I_t:= \int^t_0 \d X(s)=X(t)-X(0)
$$
when $X$ is quasi-local.
 The following bound is useful for large $t$. 
\begin{lemma} \label{lem: large time} For any $k$, 
$$
\E[ I^k_t ]^{1/k}  \leq C + Ct \big(\E[ \str f^{k}\str )^{1/k} +\sum_{x,\eta}\E [ \str g_{x,\eta}\str^k] \big)
$$
for  constants $C$ independent of $f,g$. 
\end{lemma}
Note that $f,g$ and $\tau_{Y_t} f, \tau_{Y_t} g$ have the same distribution under $\mathbf P_{\Bp, 0}$ so the RHS plays a similar role for quasi-local processes and quasi-local centered at a tagged particle.
\begin{proof}
Without loss of generality, we may assume $f, g\geq 0$.  Using the stochastic integral representation of $I$, 
$$
\d I^k_t  = \sum_{l=1}^k C(l) (\d X_t)^l  I^{k-l}_{t-}
$$
where $C(l)$ are combinatorial factors and 
\[
(\d X_t)^l:=\langle g^l, \d N_t\rangle \text{ for $l>1$.}
\]
Taking expectations, we have
\beq\label{eq: diff of moment}
\E[\d I^k_t]  = \sum_{l=1}^k C(l)   \E [( \delta_{l,1}f \d t +  \sum_{x,\eta}g_{x,\eta}^l \d t)   I^{k-l}_t].
\eeq
Applying Holder's inequality  to each term, with $1/p(l)+1/q(l)=1$ and $(k-l)q(l)=k$, we get
 $$
\frac{\d}{\d t}\E[ I^k_t]  \leq  C \E[ \str f\str^{k}]^{1/k}(1+ \E[I^{k}_t])^{1-1/k} +  C \sum_{l=1}^k  \sum_{x,\eta}  \E [ \str g_{x,\eta}\str^k]^{l/k}(1+  \E[ I^{k}_t])^{1-l/k}.
$$
Multiplying both sides by $\frac 1k (1+ \E[I^{k}_t])^{1/k-1}$ and using the fact that $(1+ \E[I^{k}_t])^{(1-l)/k} \leq 1$ leads to a differential inequality which can be integrated. The lemma follows. 
\end{proof}
For small $t$ we have a complimentary bound.
\begin{lemma}\label{lem: small time} 
$$
\E[ I^2_t] \leq  C ( z+z^2), \qquad \text{with }\,  z= t  \big( \norm f\norm_{L^2}+ \sum_{x,\eta} (\norm g_{x,\eta}\norm_{L^2} + \norm g_{x,\eta}\norm^2_{L^2}) \big)
$$
for $C$ independent of $f,g$. 
\end{lemma}
\begin{proof}
We use \eqref{eq: diff of moment} for $k=2$ and we integrate the differential inequalities in the two regimes $\E[ I^2_t] \leq 1, \E[ I^2_t] >1$, leading to the bounds $z^2,z$, respectively. 
\end{proof}

\begin{proof}[Proof of \Cref{lem: moment for use}]
The first inequality is immediate from \Cref{lem: large time}. The second follows from \Cref{lem: small time} as well by replacing $X$ by $X-X_R$ (so that  the corresponding $f,g$ are small by \Cref{ass: local}.   The fact that \Cref{lem: small time} deals only with  $L^2$-bounds is bypassed by estimating 
$$
\norm \str \nu-\nu_R \str (I)) \norm_{L^q}  \leq   \norm \str \nu-\nu_R \str (I) \norm_{L^2} \big(\norm \str \nu \str (I)) \norm^{1-1/q}_{L^{2q-2}}+\norm \str \nu_R \str (I)) \norm^{1-1/q}_{L^{2q-2}} \big)
$$
for $I=[0,t]$. The second factor is then estimated by \Cref{lem: large time}. 
\end{proof}

\section{Motion of discrepancies}
\label{S:D}
In this section, we deal throughout with processes taking values in $\Om^2$, or $\Om^2\times \Z$ when also considering tagged particles.   Pairs of spin configurations are denoted by  $\si=(\si^1,\si^2)$ with $\si^i \in \Om$.  A site $x$ where $\si^1(x)\neq\si^2(x)$ is said to host a 'discrepancy', and we say the discrepancy is of sign $+$ when $(\si^1(x)\si^2(x))=(+,-)$ and it is of sign $-$ when $(\si^1(x)\si^2(x))=(-,+)$. 
Let $D(\sigma)$  denote the position of the left most discrepancy of $\si$, i.e. 
$$
D(\sigma) :=  \inf \{x: \si^1(x)\neq \si^2(x) \}.
$$
We will always consider initial measures on $\Om^2$ so that $D>-\infty$ almost surely.
For $S \subset \Z$, let $\mu_S$ be the initial measure on $\Om^2$ defined by the following conditions:
\beq
\label{Eq-muS}
\begin{aligned}
& \text{$\si^1$ and $\si^2$ are $\Bp$ distributed},\\
& \text{For $x\in S^c$, $\sigma^1(x)=\sigma^2(x)$},\\
& \text{For $x\in S$, $\sigma^2$ is independent of $\si^1$}.
\end{aligned}
\end{equation}
That is, the measure $\mu_S$ places possibly discrepancies in all $x \in S$. 
The coupling construction defines a dynamics on discrepancies.  For example, let the configuration be $\sigma^1$(above), $\sigma^2$(below): 
$$
\begin{array}{l}   +++-+++++++-     \\  +++-+-+++++-  \\  \cdot\cdot\cdot\cdot\cdot\cdot\cdot\cdot x\,\, \,  y\,\, \,  z \cdot\cdot\cdot\cdot\cdot\cdot\cdot \cdot w    \end{array}
$$
The discrepancy (of sign $+$) sits at $y$. If the first clock ring (locally) is at $x$, then the discrepancy will move to site $w$. If the first clock ring is at $y$, then it will move to $z$ or $w$, depending on the relevant random variable $U$. Other clock rings do not move the discrepancy. In fact, it is guaranteed that a clock ring on the site of the discrepancy and a clock ring on the site left to it will move the discrepancy forward by at least one site.     In case there is more than one discrepancy around, the picture is slightly more complicated. Discrepancies of type $+$ can annihilate with discrepancies of type $-$ (they cannot cross each other) and discrepancies of the same type can possibly cross. What the latter means (to have a crossing of discrepancies) is a matter of convention. We will never need such considerations, and don't sort this out. For us, it is important to realize (by inspection of possibilities) that $1)$ the motion of an isolated discrepancy is independent of the presence of other discrepancies as long as it does not collide with or cross (or is crossed by) any of them, and 2) for the leftmost discrepancy $D$, it is in any case true that a clock ring on or left to that discrepancy will move it by at least one site.   This leads to an immediate proof of the following bound: 
\begin{proposition}[Linear Displacement of Minimal Discrepancy \textrm{I}]
\label{P:Exp}
There are constants $c, C>0$ such that for any $x \in \Z$ and all $t>0$,
\[
\mathbf P_{\mu_{[x, \infty)}}\left( D(\sigma_t)- x< {ct} \right) \leq C e^{-c t}.
\]
\end{proposition}

Whenever the tagged particle is involved, we need the following tweak of the above estimate, showing that discrepancies run away from particle with a positive relative speed.  Its proof appears in the next subsection. {We write $\mathbf{P}_{\mu_S, 0}$ for the coupled process started from the coupling measure $\mu_S$ conditioned on $\sigma^1_0(0)=+$. That is, the convention is that the tagged particle is placed in the first configuration $\sigma^1$. Therefore, we take $Y_t=Y_t(\sigma^1)$.  }
\begin{proposition}[Linear Displacement of Minimal Discrepancy \textrm{II}]
\label{P:Exptag}
There are constants $c, C>0$ such that for any {$x >0$} and all $t>0$,
\[
\mathbf P_{\mu_{[x,\infty)}, 0}\left( (D(\sigma_s)-x)- (Y_s-Y_0) < c s \,\, \text{for some $s\geq t$} \right) \leq C e^{-c t}.
\]
and $D(\sigma_t)-Y_t>0$ for all $t\geq0$, with probability $1$.
\end{proposition}

% Even though the above statements (and some upcoming ones) were proven only for the leftmost discrepancy, we point out the following:
%Let $\sigma=(\sigma^1,\sigma^2)$ be a configuration for which the dynamics can be defined and let  $ \mathbf{D}(\sigma)$ be a finite set. Then the evolution of any  discrepancy $x \in  \mathbf{D}(\sigma)$ is well-defined up to the stopping time $t_*$ where it meets another discrepancy (or annihilates with it). Moreover, up to that time, its movement is independent of the other discrepancies }

\subsection{Tagged Particles}\label{sec: tagged}
We prove here \Cref{P:Exptag}.
Let $(\sigma^1_0, \sigma^2_0)$ be two initial configurations with $\sigma^1_0(0)=+1$.  
The proof relies on the introduction of a pair of orderings associated to the particles of $\si^1$.  The first (resp. second) ordering labels the $+$ (resp. $-$) particles relative to one another.  The orderings are defined at $t=0$ and preserved in time according to the "push" dynamics.  To order the $+$-particles at $t=0$ we use the notation $i^+$ with $i \in \Z$.  We set $0^+=Y_0=0$ and label the $i$'th particle to the right or left of $0^{+}$ by $i^{+}$ depending on whether $i$ is respectively positive or negative.  
 We shall denote by $Y_t^{i^{+}}$ the position \textit{in $\Z$}  at time $t$ of the particle labeled by $i^{+}$.  An analogous ordering of the $-$-particles is fixed once we declare $0^-$ to be the first particle left of $0$ at $t=0$.

Next, we define locations in these orderings for the discrepancies appearing in $\si_0$.    Recall that a discrepancy can be either of sign $+$ or of type $-$ and its sign is conserved throughout its evolution, though, as already remarked,  opposite discrepancies can annihilate.  

Suppose there is a $\pm$-discrepancy at $x$ at $t=0$.  Assuming that it did not by the time $t>0$, denote its location \textit{in $\Z$} by $\frak d^x_t$   We'll give the another 'location' of a $+$-(resp. $-$-)discrepancy by specifying the label $i^+$ (resp. $i^-$) of the $+$-(resp. $-$-)particle the discrepancy sits on.
That is, we set $d^x_t= i^{\pm}$ where $i$ is such that $
\frak d^x_t=Y^{i^{\pm}}_t$. As long as the discrepancy is isolated, it is easy to see that $d^{x}_t$ either increases or stays constant when a clock ring affects the discrepancy.  In fact, if the clock at $\frak d^x_t$ rings and the relevant $U$-variable dictates the $\mp$ particle to move, then $d^{x}_t$ is guaranteed to increase by at least one.  This means that the increase of $d^x_t$ may be stochastically bounded from below by a rate $\min(\lambda_+, \lambda_-)$ Poisson process.  

There is ambiguity in this reasoning when other discrepancies touches are present unless we focus on (only) the leftmost discrepancy of type $\pm$.
In that case, among all potential outcomes, the only one requiring further explanation is when the leftmost $\pm$-discrepancy annihilates  with one of opposite type.  In that, case one of the discrepancies to its right becomes the leftmost discrepancy (or it is assigned the value $\infty$, if there is no other discrepancy of the same type).
The foregoing discussion, with $D^\pm(\sigma_t)$ denoting the position, \textit{in the $\pm$ ordering}, of the leftmost discrepancy of type $\pm$, proves the following: 
\begin{lemma} \label{lem: doutrunstag}
There exist $C,c>0$ such that
\label{L:Disc1}
\[
\mathbf P_{\mu_{[x,\infty)}}\left(D^{\pm}(\sigma_t)-D^{\pm}(\si_0)< c t\right) \leq C e^{-c t}.
\]
\end{lemma}
We can now proceed with
\begin{proof}[Proof of   \Cref{P:Exptag}]
As long as the left-most discrepancy $D$ has sign $+$, the claim is easy:  The position of the discrepancy in the $+$ ordering is linearly increasing by \Cref{lem: doutrunstag}, whereas the position of the tagged particle in the $+$ ordering is constant. This also implies a linearly growing distance on the lattice. 
When the left-most discrepancy has sign $-$, it takes valued in a different ordering than the tagged particle, so the above argument fails. However, since the tagged particle is to the left of all discrepancies,  its motion in $\sigma^1$ and $\sigma^2$ is the same. Therefore, one may now reverse the roles of $\sigma^1$ and $\sigma^2$, thus flipping the sign of the discrepancy so it takes values in the same ordering as the tagged particle. 
\end{proof}
.

%\begin{proof}
%We observe that, when finite, $d^x_t$ only takes values in one of the two linear orders for all $t$.  Which of the orderings it does take values in depends only on whether $\sigma^1_x(0)= \pm$.  As a consequence of the dynamics, $d^x_t$ and hence $D^{\pm}(\si_t)$ is increasing.  Further, $D^\pm(\sigma_t)$ is guaranteed to increase under the circumstance that the clock its location in $\Z$ and that $\mp$-particles at that vertex move in the coupling dynamics.  Our claim now follows since the displacement of $D^\pm(\sigma_t)$ in the $\pm$ ordering may be stochastically bounded from below by a rate $\min(\lambda_+, \lambda_-)$ Poisson process.
%\end{proof}

\subsection{Upper Bound on Speed of Discrepancies}
Above, we have argued that discrepancies move at least linearly to the right/away from tagged particles. Now we provide upper bounds. 
\begin{lemma}  \label{lem: max v dis}
For any $t \geq 1, R\geq 0$, we have 
$$
\prob_{{\mu_{\{x\}}}} (\si^1_0(x)\neq \si^2_0(x), \: D(\si_t)-x  \geq R)  \leq  C \e^{-c (R/t)^{1/2}}
$$
\end{lemma}
The same reasoning can be used to prove bounds on the displacement of the tagged particle. 
\begin{lemma}\label{lem: max v tag}
For any $t \geq 1, R \geq 0$, we have 
$$
\prob_{\Bp, 0} (\str Y_s-Y_0 \str \geq R \,\, \text{for some $s\geq t$})  \leq  C \e^{-c (R/t)^{1/2}}
$$
\end{lemma}
To prove these results we need an \textit{a-priori} flux bound:
Let us define the counting processes
\beq \label{eq: instant current}
J_{x}(t):=\sum_\eta  \sum_{y<x} \int_0^t   \chi^{-\eta}_{[y,x-1]} (\si_{s-})  \chi^{\eta}_{y} (\si_{s-})  \d N_{y,\eta}(s)
\eeq
This process records the total number of particles which jump from $(-\infty, x)$ to $[x, \infty)$ in the time interval $[0, t]$.  The following bound was proved in \cite{crawford2015toom}, see Lemma 4.7.
\begin{lemma}[\textit{A-Priori} Flux Bound]
\label{L:Flux}
There are constants $C, \gamma>0$, depending only on $\lambda_{\pm}, p$, such that
\[
\E_{\Bp}[\e^{\gamma J_{x}(t)/t}] < C, \qquad \text{for any $t>0$}
\]
\end{lemma}

\begin{proof}[Proofs of \Cref{lem: max v dis,,lem: max v tag}]
For concreteness, we restrict attention to the proof of \Cref{lem: max v dis}, the argument being similar in the remaining case.   Let us denote by $l_{x}$ the number of spins to the left of $x$ (including $x$) of the same sign as $\sigma(x)$.  Similarly, $r_x$ is the number of like spins to the right, starting at $x$. By definition $l_x,r_x \geq 1$.

There are two ways a discrepancy at $x$ can move:  The first way is that the exponential rate one clock on the vertex it occupies rings. 
The other way is if one of the exponential rate one clocks at $x-l_x+1,\ldots, x-1 $ rings.  Hence the local rate of moves is bounded by $l_{D_t}$. If such a move occurs, the jump length is bounded by $r_{D_t}$. 
%  In both cases the displacement is bounded by $R_{D_t}+1$, but may be smaller. Thus, the local rate of movement is bounded above by $(|L_{D_t}|+1)(R_{D_t}+1)$. 
   So, if we can bound the size of stretches of like spins that the discrepancy encounters, we can bound the speed of the discrepancy.  

Let $ E$ be the event that a stretch of spins of length at least $L$ occurs in the spatial interval $[x-R, x+R ]$ in the time interval $[0, t]$, hence not necessarily only adjacent to $D_s$ for $s \in [0, t]$.  
We will nevertheless find a good bound on $E$ and then estimate the motion on $E^c$ straightforwardly.
The parameter $L$ will be fixed, depending on $R,t$ at the end of the proof.  We divide $[0, t]$ into $t$ intervals of length of order $1$ and divide $[x-R, x+ R ]$ into blocks of length $L$. We enumerate the corresponding spacetime rectangles of $[0, t] \times [-R, R ]$ by $(B_{j})_{j=1}^J$, where $J= O(R t)$. We write $B_j= [s_j, t_{j}) \times [a_j, b_j)$ and consider the events
\[
E_{j}=\{ \text{$[a_j, b_j)$ has a stretch of $L/2$ like spins at some $t\in [s_j, t_{j})$}\}.
\]
\[
F_{j}=\{ \text{$[a_j, b_j)$ has a stretch of $L/4$ like spins at $s_j$}\}.
\]
In order for $E_j$ to occur, either there must already be a stretch of length $L/4$ present at time $s_j$, i.e.\ $F_j$ occurs, or at least $L/4$ particles must cross some vertex $x \in [a_j, b_j)$ in the (small) time interval $[s_j, t_{j})$.  Both of these possibilities are unlikely: 
A large deviation estimate for $\Bp$ yields
\[
\prob(F_j)\leq C\e^{-cL}.
\]
and the flux bound \Cref{L:Flux} bounds the probability that $L/4$ particles crossed a vertex, i.e.\
\[
\prob(E_j\str F^c_j)\leq C\e^{-cL}.
\]
Hence we conclude that $\prob(E_j) \leq C\e^{-cL}$  and hence $\prob(E) \leq C(tR)\e^{-cL}$. It remains to estimate the speed of the discrepancy condition on $E^c$. As explained above, the distance traveled is now bounded above by $LN^{(L)}_{t}$ with $N^{(L)}_t$ a Poisson process with intensity $L$. 
Large deviation estimates yield that $\prob(LN^{(L)}_{t} \geq R) \leq C\e^{-c R/L} $ provided that $R  \geq C t L^2$. 
Collecting the estimates, we obtain
\[
\prob_{{\mu_{\{x\}}}} (\si^1_0(x)\neq \si^2_0(x), \: D(\si_t)-x  \geq R) \leq   C(tR) \e^{-cL}  + C \e^{-cR/L}, \qquad \text{for $R \geq Ct L^2$}
\]
which is optimized to give a bound $C \e^{-c\sqrt{R/t}} $, provided $t \geq 1$.
\end{proof}

\subsection{Decay of Correlations}

%Let $S_i \subset \R^d$ be some compact Borel sets and recall that $\caM(S_i)$ is the Banach space of finite signed measures with the TV norm.  
For a random measure $\mu=\mu(\si)$ we write 
$\tau_{-x}\mu(\sigma)= \mu(\tau_x\sigma)$, i.e.\ same convention as for number-valued random variables. {Also, as a natural extension of our previous notation we will say $\mu\in \caB(\sigma(x): x\in A)$ if for every $f\in C_b(S_i)$, the variable $\int f \textrm{d} \mu\in \caB(\sigma(x): x\in A)$.}

\begin{lemma}[Exponential Decay of Correlations]\label{lem: decay of correlations}
Let  $U, V $ be random measures on compacts $S_1,S_2$, respectively  such that $U \in \caB(\sigma(x): \str x \str  \leq M)$ and $V\in \caB(\sigma(x):  x   \leq M)$.    Then \begin{multline}
\left\str \E_{\Bp,0}[\tau_{Y_t}V(\sigma_t) U(\sigma_0)] -     \E_{\Bp,0}[V(\sigma_0)]  \E_{\Bp,0}[ U(\sigma_0)] \right\str(S_1\times S_2) \\[1mm]
\leq   C \norm \str U\str(S_1)  \norm_{L^4}  \, \norm   \str V\str(S_2)  \norm_{L^4} \,\e^{-c(t-M)}
\end{multline}
where all $\str \cdot \str(S)$ stand for the variation on $S$.   The same bound holds if we replace $\tau_{-Y_t}V$ with $V$ and $\E_{\Bp,0}$ with $\E_{\Bp}$ (i.e.\ the case with no tagged particle).
\end{lemma}
 The statement in the absence of a tagged particle is simpler to prove. In fact, a weaker version applying to functions rather than measures, appears already in \cite{crawford2015toom}. Thus, we explicitly prove here only the decay of correlations in the presence of a tagged particle.  There are some technical complications, mostly due to the fact that if one tries to couple two tagged particles in two different environments, they will not necessarily lie on the same vertex in $\Z$ after all discrepancies move to the right of them. To circumvent this difficulty, the idea is to focus on a tagged particle that starts to the left of all discrepancies.
 
 \begin{proof}[Proof of \Cref{lem: decay of correlations}]
Without loss, let $V$ to be of zero mean. 
 Let 
$$
-M'(\sigma_0)=\inf \{x < -M: \sigma_0(x)=+ \},
$$
that is, $-M'$ is the position of the rightmost $+$-particle to the left of $-M$.  We set
\[
 Z_t(\si)=Y^{-M'(\sigma_0)}_t,
\]
so that $Z_t(\si)$ is position at time $t$ of the tagged particle started from $-M'$.

Given $V:\Omega \rightarrow \caM(S_1)$ and $n \in \bbN$, let  $V(\sigma, n)$ be the measure $V$ shifted to the $n$'th $+$-particle right of the origin, i.e.
$$
V(\sigma, n):= \tau_{-\tilde n} V(\sigma), \qquad \text{with $\tilde n=  \min \left\{ m: \sum_{i=1}^{m}\chi(\sigma(i)=+)=n \right\} $}
$$
Then we have the identity
$$  
 \tau_{-Y_t}V(\sigma_t)    =    \tau_{-Z_t(\sigma)}V(\sigma_t,N(\si_0)).
$$
where $N(\si_0)$ is the number of +-particles between $-M$ and $0$.   Both $M'$ and $N$ are random and depend on $\sigma_0$. Crucially however, they are independent of one another under the measure $\Bp^0:=\Bp(\,\cdot \,  \str  \sigma_0(0)=+)$.
 
Let us consider the coupling measure $\mathbf P_{\mu[-M,\infty),0}$, as defined at \cref{Eq-muS} except that $\sigma^1_0$ is conditioned to have  $\sigma^1_0(0)=+$. 
Note that $Z_t(\si^1)=Z_t(\si^2)$ because $M'(\sigma_0^1)=M'(\sigma_0^2) $ and a tagged push particle started to the left of all discrepancies can never catch up with the discrepancies, see \Cref{P:Exptag}.
Let $A$ be the event that at time $t$ the leftmost discrepancy is to the right of $Y_t(\sigma^1)+M$, where $Y_t(\sigma^1)=Y^0_t(\sigma^1)$ is the tagged particle \emph{started from the origin}.  On $A$, we have
\beq \label{eq: def h}
 \tau_{-Z_t(\si^1)}V(\sigma^1_t,N(\sigma^1_0))=    \tau_{-Z_t(\si^2)}V(\sigma^2_t,N(\sigma^1_0))=:h(\sigma^2,\sigma^1_0)
\eeq
since  $Z_t(\si^1)=Z_t(\si^2)$ and $\si^1_t(x)=\si^2_t(x)$ for $x$ smaller than the leftmost discrepancy. 

Recall that by \Cref{P:Exptag}, the event $A$ occurs with probability at least $1-C\e^{c(M-t)}$.  
Using \Cref{eq: def h} and $1=1_A+1_{A^c}$, we get
\begin{multline}  \label{eq: split fg}
\E_{\Bp,0}[\tau_{-Y_t}V(\sigma_t) U(\sigma_0) ] - \E_{\mu[-M,\infty),0}\left[  h(\sigma^2,\sigma^1_0) U(\sigma_0^1) \right]   \\[1mm]
=     -   \E_{\mu[-M,\infty),0}[1_{A_c} \left(  \tau_{-Y_t}V(\sigma^1_t) U(\sigma_0^1) -h(\sigma^2,\sigma^1_0) U(\sigma_0^1) \right) ]
\end{multline}
The second term on the left hand side may be re-expressed as 
$$ 
\int \d \Bp^0(\sigma^1_0) U(\sigma_0^1)   \E_{\mu[-M,\infty),0} \Big[   h(\sigma^2,\sigma^1_0) \big\str \sigma_0^1(x), \str x \str \leq M   \Big] 
$$ 
The random variable $ h(\sigma^2,\sigma^1_0)$ depends on $\sigma_0^1$ only through $N$, so we can conclude that
\begin{multline}
 \E_{\Bp,0} [    h(\sigma^2,\sigma^1_0)  ] =   \sum_{n\in \N} \mathbf \chi\{N(\si_0)=n\} \E_{\Bp,0} [   \tau_{-Y_t}V(\sigma_t,n)  ] \\
 =   \sum_{n\in \N}  \E_{\Bp,0} [   \tau_{-Y_t}V(\sigma_t)  ]   =0.
\end{multline}
The second equality follows from translation invariance and the third follows since $V$ is of zero mean.  It follows that the second term on the left hand side in \Cref{eq: split fg} vanishes and to conclude the proof, we need to estimate the total variation of the right hand side in \Cref{eq: split fg}, which is of the form $\widetilde\E(J 1_{A^c})$ with $J$ a measure and $\widetilde \E=\E_{\mu[-M,\infty),0}  $.  We use 
$$ \str \widetilde\E(J 1_{A^c}) \str \leq  \widetilde\E( \str J \str  1_{A^c})  \leq \widetilde\E(\str J \str^2)^{1/2} (\widetilde{\mathbf{P}}(A^c)^{1/2}, $$
with $\str \cdot \str$ denoting total variation and $\widetilde \E=\E_{\mu[-M,\infty),0}  $  As already remarked, the probability of $A^c$ is exponentially small, so we just need to bound $\widetilde\E(\str J \str^2)^{1/2}$, which goes as follows:
$$
\widetilde\E(\str J \str^2(S_1\times S_2))^{1/2} \leq  2   \E_{\Bp,0}(\str U \str^4(S_1))^{1/4}   \E_{\Bp,0}(\str V \str^4(S_2))^{1/4}
$$
where we used stationarity of the process seen from the tagged particle.
\end{proof}

\section{Time-Reversal and the Adjoint Process}

\subsection{The Time-Reversal Map}
%  In particular, for any times $t>t_i>0$, functions $f_i \in C(\Om)$ and for $\si$ $\Bp$ a.s., 
%\beq
%\label{AdjCondExp}
%   \E\left[\prod_{i=1}^{l}f_i(\sigma_{t-t_i}) \big\str \sigma_t=\si\right]  = \E^{*}_{\si}\left[{ \prod_{i=1}^{l}f_i(\sigma_{t_i})} \right]
%\eeq
Let us fix some time $\tau>0$ and define the time-reversal map 
$\sigma \mapsto \tilde \sigma $  from $D_{\Om}([0,\tau])$ to $D_{\Om}([0,\tau])$ by 
$$
\tilde \sigma_{s}:=  \sigma_{(\tau-s)_-}, \qquad 0 \leq s \leq \tau
$$
This map is measurable  and is one-to-one on a set of full $\prob_{\Bp}$ measure.
Let $F  \in \caB(D_\Omega([0,\tau]))$ and consider the lift of the time-reversal map to functions $F \mapsto \tilde F$:
$$
\widetilde F(\sigma) :=  F(\tilde\sigma).
$$
For each of our quasi-local processes $X_t$, we now have a time-reversed process $\tilde X_t$ satisfying
$$
\tilde X_t(\sigma)-\tilde X_s(\sigma) :=  -( X_{\tau-s}(\tilde\sigma)- X_{\tau-t}(\tilde\sigma)), \qquad 0\leq s \leq t \leq \tau,  
$$
{It is instructive to take $X_t=\int_0^t~\chi\{\si_{s-}(x)=\eta\} \textrm{d} N_{x,\eta}(s)$.}  In this case, comparing $\sigma, \tilde\sigma$ at jump times, we get
$$
\d \tilde X=  - \sum_{y>x}  \chi^{-\eta}_{x} \chi^{\eta}_{(x,y]}(\si)  \d N_{y,\eta}(t).
$$
This allows us to deduce that the mapping $X \to \tilde X$ maps quasi-local processes into quasi-local processes. 
The thing to keep in mind is that an arrival of $N_{x,\eta}$ at time $s$ causing a jump for the process $X$ corresponds to an arrival of $N_{y,\eta}$ at time $t-s$ for the process $\tilde X$ where $y= \min (z: z>x, \sigma_{s-}(z+1)=-\eta)$.  
More generally, with $X$ determined by $(f,g_{x,\eta})$, the map $X \mapsto \widetilde X$ corresponds  to the map $(f,g_{x,\eta}) \to (\tilde f,\tilde g_{x,\eta})$ with
$$
\tilde f=-f, \qquad  \tilde g_{y,\eta} = \sum_{x<y}  g_{x,\eta}  \chi^{-\eta}_{x} \chi^{\eta}_{(x,y]}.
$$
The data $(\tilde f,\tilde g_{x,\eta})$ satisfy all necessary requirements: 
\begin{lemma}\label{lem: loc of tilde}
If, as assumed throughout, $(f,g_{x,\eta})$ are such that all three Assumptions â\ref{ass: norms} \ref{ass: local} and \ref{ass: real} are satisfied, then they are satisfied as well for $(\tilde f,\tilde g_{x,\eta})$. 
\end{lemma}
The straightforward verification of this lemma proceeds by using Holder inequalities and the fact that $\norm \chi^{\eta}_{(x-r,x]} \norm_{L^q(\Bp)} \leq C \e^{-c r}$ for any $q >0$.

\subsection{The Adjoint Process $\bf{P}^*_{\Bp}$}
\label{S:Adjoint}

Let us denote by $\E^*_{\Bp}$  the expectation started from $\Bp$ of a left-moving Toom interface.  Thus when the $N_{x,\eta}$ clock rings and $\sigma(x)=\eta$, we exchange the values of $\sigma(x),\sigma(y)$ with $y:= \max(z<x: \sigma(z)\neq \sigma(x))$.    The left-moving process started from $\Bp$ is constructed analogously to the right moving process  and again $\Bp$ is an invariant measure.  The process can be started from $\si$ $\Bp$-almost surely, and we denote its expectation by $\E^*_\si$.  
The relation to the time-reversal map introduced above is that
\beq
\label{AdjCondExp}
   \E_{\Bp}[ F ]  = \E^{*}_{\Bp}[ \tilde F ] 
\eeq

Let us briefly sketch the verification of \Cref{AdjCondExp}. First, using \Cref{R:Adj} one verifies \Cref{AdjCondExp} for functions on $D([0, \eps])$ (note that on a finite cycle the corresponding statement is direct).  Then using the Markov property and induction, one extends to functions on $D([0, \tau])$ for arbitrary $\tau$

Here is the induction step:
Assume \cref{AdjCondExp} for functions of $D([0, t])$.  We extend it to functions of $D([0, 2t])$.  Let $s\in [0, 2t]$ and let $f, g$ be bounded measurable functions.
Then $$\E[g(\si_0)f(\si_{s})]=\E[g(\si_0)\E_{\si_{s/2}}[f(\si_{s/2})]]=\E^*[g(\si_{s/2})   \E_{\si_{0}}[f(\si_{s/2})]].$$  The first equality  follows from the Markov property while the second follows from the induction hypothesis for \cref{AdjCondExp}.  Note here that the outer expectation corresponds to the left moving process while $\E_{\si_{0}}[f(\si_{\eps})]$ corresponds to the right moving process.  
Using the Markov property again (for left moving process)
the RHS is $\Bp(\si_0)[E^*_{\si_{0}}g(\si_{s/2})   E_{\si_{0}}[f(\si_{s/2})]]$.  We are then done by symmetry ($E^*[g(\si_{2\eps})f(\si_{0}]$ yields the same expression).   The argument for a general finite product at different times in $[0, 2t]$ is similar.  Then we conclude the induction step by density argument (or by the Monotone Class Theorem).

If we want to include the tagged particle, we begin by considering functions $f_i$ on the extended state space $\Om \times \Z$.  It simplifies matters to assume that each $f_i$ is translation covariant, i.e.\ $f_i(\sigma,y)=f_i(\tau_x\sigma,y-x )$, in which case \Cref{AdjCondExp} is upgraded to
\beq
\label{AdjCondExptagged}
   \E_{\Bp,0}[ F ]  = \E^{*}_{\Bp,0}[ \tilde F ] 
\eeq

Let us fix a time $s$ and we consider  two $L^2$ functions $F_1,F_2$ where $F_1\in \caB(\sigma_t, t \in [0,s))$ and $F_2\in \caB(\sigma_t, t \geq 0)$.
Let
$$
G(\sigma) :=  \E^*_{\sigma}( \widetilde F_1 ).
$$
Note that the $\widetilde{\cdot}$ operation depends on a fiducial point $\tau$, which is understood here to be $\tau=s$. 
\begin{lemma}\label{lem: reversal}
With $s, F_1, F_2, G$ as above
$$
\E[F_1\:  F_2\circ \theta_s] =   \E[G F_2 ].
$$ 
\end{lemma}
\begin{proof}
We have
\begin{eqnarray}
\E_{\Bp}[F_1\:  F_2\circ \theta_s] &=& \E_{\Bp}[F_1\: \E_{\si_s}[F_2]]\\
&=&   \int \d\Bp(\si)   \E_{\Bp}[F_1\str \sigma_s= \si]   \E_{\si}[ F_2]  \\
&=&     \E_{\Bp}\left[ G  F_2 \right]
\end{eqnarray}
Here the first  equality is due to the Markov property, the second is due to stationarity of $\Bp$ and the definition  of the conditional expectation.  The third equality follows from the fact that  $\E^*_{\sigma}( \widetilde F_1 )$ is a version of   $\E_{\Bp}(F_1 \str \sigma_s=\sigma)$ cf. 
\Cref{AdjCondExp}.
\end{proof}

\subsection{An Application}\label{S:UBApp}
 
To foreshadow future applications, we use  \Cref{lem: reversal} to obtain identities between measures generated by $\ell$ quasi-local processes $X^{(j)}(t)$. Let us first assume that $(X^{(j)})_{j=1,\ldots, \ell}$ are quasi-local w.r.t. to the origin. Let $l \in [1, \ell]$ and observe that
$$
\E_{\Bp}[   \d X^{(1)}(t_1)\ldots   \d X^{(\ell)}(t_\ell) \big ] =  \textrm{d} t_{l}\E_{\Bp}[H(\d t_{1,\ldots,l-1})   \d X^{(l+1)}(t_{l+1})\ldots   \d X^{(\ell)}(t_{\ell})]
$$
where $H(\d t_{1,\ldots,l-1})=H(\sigma, \d t_{1,\ldots,l-1})  $ is the random variable on $\Om$, defined $\Bp$-a.s., with values in measures on $\caM(\caS_{l-1})$, cf.\ \Cref{lem: def measures} given by 
$$
H(\sigma, \d t_{1,\ldots,l-1})=  \E_{\Bp}[   \d X^{(1)}(t_1)\ldots  \d X^{(l-1)}(t_{l-1})\,\str\,  \sigma_{t_l}=\sigma]
$$

It is useful  to rewrite this formula  using the adjoint process. 
Consider the change of variables
$$
s= t_{l}, \qquad    w_{j}= t_{l}-t_{l-j}, \quad \text{for $j =1,\ldots, l-1$} \quad  u'_{j}= t_{j+l}-t_{l}, \quad \text{for $j =l+1,\ldots, \ell$}.
$$
We will work in \Cref{sec:-app-to-h,S:CLT} with this change of variables.

Then using stationarity (and abusing the notation for $H$),  on the set $\caA:=\{ 0< w_1< \dotsc < w_{l-1}<s, 0< u'_1< \dotsc u'_{\ell-l}\}$
\beq
\label{E:COV}
\E_{\Bp}[   \d X^{(1)}(t_1)\ldots   \d X^{(\ell)}(t_\ell) \big ] \mapsto  \textrm{d} s\E_{\Bp}[H(\si, \d w_{1,\ldots,l-1})   \d X^{(l+1)}(u'_1)\ldots   \d X^{(\ell)}(u'_{\ell-l})].
\eeq
Here
the measure $H(\sigma, \cdot)$ satisfies
\begin{align}
\label{E:H}
H(\sigma, \d w_{1\ldots l-1}) =   &\tilde f(\sigma) \E^*_{\sigma}(\d \tilde X^{(l-1)}(w_1)\ldots   \d \tilde X^{(1)}(w_{l-1})   ) \\
& +   \sum_{x,\eta, r}  \chi^{-\eta}_{r} \chi^\eta_{(x-r,x]} \tilde g^{(l)}_{x,\eta}(\sigma) \E^*_{\sigma^{x-r,x}}(\d \tilde X^{(l-1)}(w_1)\ldots   \d \tilde X^{(1)}(w_{l-1})) \nonumber
\end{align}
where $r$ ranges over $r=0,1,2,\ldots$ and $\sigma^{x-r,x}$ is obtained from $\sigma$ by exchanging $\sigma(x-r),\sigma(x) $.

Eventually, we are interested in the total variation over the set $E_l(T, \kappa)$ (see \Cref{E:EE}), which  under the change of variables gives the restrictions $w_{l-1}<T$ and $u'_1>T^\kappa$ on $\caA$ (among other conditions).  It is convenient on $\{u'_1>T^\kappa\}$ to change variables once more.  Applying the Markov property and setting $u_i= u_i'-T^\kappa$, our measure transforms into
\beq
\label{E:COV2}
\E_{\Bp}[   \d X^{(1)}(t_1)\ldots   \d X^{(\ell)}(t_\ell) \big ] \mapsto  \textrm{d} s\E_{\Bp}[H(\si, \d w_{1,\ldots,l-1}) K(\sigma_{T^\kappa},\d u_{1,\ldots,\ell-l } ) ]  
\eeq
where  the measure $K$ is given by 
$$
K(\sigma, \d u_{1,\ldots,\ell-l })  = \E_{\sigma}[  \d X^{(l+1)}(u_1) \ldots   \d X^{(\ell)}(u_{\ell-l})]
$$
and 
\beq
\label{E:EEE}
E_l(T, \kappa)\mapsto \{ 0< w_1< \dotsc < w_{l-1}<T<s, \: 0< u_1< \dotsc u_{\ell-l}< L-T^\kappa- s\}=:E'.
\eeq
If we consider $X^{(j)}$ quasi-local w.r.t. to the tagged particle, then we can write the same formulas as above provided we replace
$\E_{\Bp}, \E^*_{\sigma}, \E^*_{\sigma^{x-r,x}}$ with $\E_{\Bp,0}, \E^*_{\sigma,0}, \E^*_{\sigma^{x-r,x},0}$ and $K$ by its natural analog depending on $\si_{T^\kappa}$ and $Y_{T^\kappa}$.

\section{Bounds on Localization of Conditional Expectations}
\label{S:LocalCE}

\subsection{A General Principle}
Let  the function $F$ satisfy $F\in \caB(\si_s(x), (x,s) \in A )$ for some Borel set $A \subset \Z \times \R_+$.  Let
$$
G(\sigma):= \E_{\sigma}(F)
$$ 
and define its local approximations by
\[
\texttt{P}_R G(\sigma): = {\Bp}\bigg[G(\cdot)|\sigma_0(j)=\sigma(j): j \in [-R, R]\bigg].
\]
We define the event 
$$
E_A:= \{ \text{$x \in \mathbf{D}_t$ for some  $ (x,t) \in A$} \}.
$$
i.e.\ the event that there is a discrepancy in $A$.  
Let $\nu_k$ be defined as the measure on pairs of spin configurations $\si=(\si^1, \si^2)$ such that:
\begin{align*}
& i) \text{ The marginal distributions of $\sigma^1$ and $\sigma^{2}$ are $\Bp$. }\\
&ii) \text{ For $j \neq \pm (k+1)$, $\sigma^1(j)=\sigma^2(j)$ $\nu_k$-a.s.}\\
& iii) \text{ For $j=\pm (k+1)$, $\sigma^1(j), \sigma^2(j)$ are independent.}
\end{align*}
We can  now state the main result of this section.
\begin{lemma} \label{L:Local2}
For any $R>0$
\[
\Bp\left[ \left(G-\texttt{P}_RG\right)^2\right] \leq C\sqrt{\E_{\Bp}[F^4]}\sum_{k \geq R} \sqrt{\mathbf{P}_{\nu_k}[E_A] }
\]
For the  tagged particle case, we get the same lemma with the replacements 
\begin{align*}
&F\in \caB(\si_s(x): (x,s) \in A ) \to F\in \caB(\si_s(x), Y_s: (x,s) \in A ),\\
&G(\sigma)= \E_{\sigma}(F)\to G(\sigma)= \E_{\sigma, 0}(F),\\
&\mathbf{P}_{\Bp}\to \mathbf{P}_{\Bp,0}. 
\end{align*}
\end{lemma}
\begin{proof}
Using the natural decomposition of in terms of martingale increments, 
$$
\Bp\left[(G-\texttt{P}_RG)^2\right] = \sum_{k \geq R}   \Bp\left[(\texttt{P}_k G-\texttt{P}_{k+1} G)^2\right]
$$
We represent 
\beq \label{eq: diff}
\texttt{P}_{k} G(\sigma) - \texttt{P}_{k+1} G(\sigma)  = {\nu_{\sigma, k}}  [  G(\sigma^1) -  G(\sigma^2)  ] = \E_{\nu_{\sigma, k}}  [  F(\sigma^1) -  F(\sigma^2)  ] 
\eeq
where, for $\si \in \Om$, $\nu_{\sigma, k}$ is the probability measure on $\Om^2$ (with configurations
$(\si^1, \si^2)$) such that:
\begin{align*}
& i) \text{ $\sigma^1(x)=\sigma^{2}(x)=\sigma(x)$ for $\str x \str \leq k$ and for $x= \pm (k+1)$, $\sigma^2(x)=\sigma(x)$.}\\
& ii)\text{ The variables  $\sigma^1(x), \sigma^2(x)$ for  $\str x \str > k+1$ and $\sigma^1(x)$ for $x=\pm (k+1)$}\\
& \quad \quad\text{ are i.i.d. They are $+1$ with prob $p$ and $-1$ with prob $1-p$.}
\end{align*}

By \Cref{eq: diff} and the Cauchy-Schwartz inequality,
$$
\Bp [(\texttt{P}_{k} G- {\tt P}_{k+1} G)^2 ]  \leq  \E_{\nu_k}[(F(\sigma^1) - F(\sigma^2) )^2].
$$
The utility of this last inequality is to reduce the derivation of error bounds in the localization of $G$ to controlling the behavior of a pair of discrepancies. 
By assumption on $F$, for any fixed $k$, 
\beq
\E_{\nu_k}[(F(\sigma^1) - F(\sigma^2) )^2]=   \E_{\nu_k}[  (F(\sigma^1)-F(\sigma^2))^2 1_{E_A)}] 
%& \leq& (\E_k( ((F(\sigma^1)-F(\sigma^2))^4 )^{1/2}  \sqrt{\mathbf{P}_{\nu_k}(E_A)}  \\
\eeq
so that 
\beq
\Bp [(\texttt{P}_{k} G- {\tt P}_{k+1} G)^2 ] \leq C\sqrt{\E_{\Bp}[F^4]} \sqrt{\mathbf{P}_{\nu_k}(E_A)}
\eeq
since the marginals of $\nu_k$ are $\Bp$.
\end{proof}

\subsection{An Application to the Random Measure $H$}\label{sec:-app-to-h}

Continuing the discussion from \Cref{S:UBApp}, we want to apply \Cref{L:Local2} to argue that the random measure $H(\sigma,\d t_{1,\ldots,l})$ can be localized in $\sigma$.  With a view toward the justification of \Cref{L:Factor1}, we shall bound the variation of $H$ on
 $$\caS_{l-1}(T) $$
 
To state the main point of this section, let us first extend the action of the projection/conditional expectation $\texttt{P}_R$ to random measures as follows.  For $\nu$ a random measure taking values in $ \caM(\caS_{\ell}(L) )$,  we define $\texttt{P}_R\nu$ by
\[
\int h \texttt{P}_R(\nu):= \texttt{P}_R\left(\int h \nu\right)
\] 
for all bounded continuous functions $h$ on $\caS_{l}(L)$.
To check that this is a meaningful definition, note that
$$
\left\str \texttt{P}_R\left(\int h \nu\right) \right\str \leq  \texttt{P}_R(\str\nu\str(\caS_{\ell}(L)))   \norm h \norm_{\infty}
$$
which is finite almost surely, since $\str\nu\str(\caS_{\ell}(L))$ is finite almost surely.  That the random measure is well-defined then follows from Lemma \ref{lem: def measures}.

\begin{lemma} \label{lem: loc of h}
Let $\kappa>1$. For $R \geq T^{\kappa}$, 
$$
\norm \str H-\texttt{P}_R H\str (\caS_{l-1}(T)) \norm_{L^q}  \leq C\e^{-cT^{\frac{\ka-1}{2}}}  
$$
\end{lemma}

Starting from the expression for $H$,  \Cref{E:H}, the main technical issue in proving \Cref{lem: loc of h} is to deal with the non-locality of the measure
$$
Z(\d w_{1  \ldots l-1}) = Z(\sigma, \d w_{1  \ldots l-1})  :=  \E^*_{\sigma}[\d \tilde X^{(l-1)}(w_1)\ldots   \d \tilde X^{(1)}(w_{l-1})   ].
$$
As such, we first study this expression separately.  
The main application of \Cref{L:Local2} is the following.
\begin{lemma}\label{lem: loc of z}
With $T,R, \alpha$ as in \Cref{lem: loc of h}, 
$$
\norm \str Z-\texttt{P}_R Z\str(\caS_{l-1}(T)) \norm_{L^q}  \leq C\e^{-cT^{\frac{\ka-1}{2}}}.
$$  
\end{lemma} 

\begin{proof}[Proof of \Cref{lem: loc of z}] 
Let us denote by $Z_{R/2}$ the result of replacing all $\tilde{X}$ by their localized versions $\tilde{X}_{R/2}$, see \Cref{S:DofC}.  Abbreviating $\caS_{l-1}(T))= \caS$ and using the triangle inequality,
\[
\norm \str Z-\texttt{P}_R Z\str(\caS) \norm_{L^q} \leq \norm \str Z_{R/2}-\texttt{P}_RZ_{R/2}\str(\caS)\norm_{L^q} + \norm \str Z-Z_{R/2}\str(\caS)\norm_{L^q} +   \norm \str \texttt{P}_RZ_{R/2}-\texttt{P}_RZ \str(\caS)\norm_{L^q}.
\]
We observe that the last contribution on the RHS is bounded by the second contribution by applying Jensen's inequality for conditional expectations.
To bound the second contribution, we have, writing $\nu^{(j)}$ for the measure $\d \tilde X^{(j)}$,  
\[
\norm \str Z_{R/2}-Z\str(\caS)\norm_{L^q}\leq C\Bigg \norm \E^{*}_{\sigma}\bigg[ \sum_{i} \str\nu^{(i)}-\nu^{(i)}_{R/2}\str[0,T]  \prod_{j \neq i} \str\nu^{(j)}\str[0,T] \bigg] \Bigg\norm_{L^q}.
\]
Using Jensen's inequality and then H\"{o}lder's inequality,  
\[
\norm \str Z_{R/2}-Z\str(\caS)\norm_{L^q}\leq C\max_{j \leq l-1} \norm \str\nu^{(j)}-\nu^{(j)}_{R/2}\str[0,T]\norm_{L^{q(l-1)}} \prod_{j=1}^{l-1} (1+\norm  \str\nu^{(j)}\str[0,T]\norm_{L^{q(l-1)}}).
\]
By \Cref{lem: moment for use}, the RHS is bounded by $C \e^{-cR}$ (which is sufficient for the claimed bound of Lemma \ref{lem: loc of z}).
Therefore we reduce the proof to providing bounds on $Z_{R/2}-\texttt{P}_R Z_{R/2}$.  

To put ourselves in the framework of \Cref{L:Local2} observe that the conclusion of said lemma applies equally to the adjoint process by symmetry and  recall that the variation of a measure on $\caS_{l-1}(T)$ may be characterized by
$$\str Z_{R/2}-\texttt{P}_RZ_{R/2}\str(\caS_{l-1}(T)) =\sup_{h: \norm h \norm_{\infty}=1} \str G_h-\texttt{P}_RG_h\str. $$
where
$G_h:= \int_{\caS_{l-1}(T)} h Z_{R/2} $ for a bounded continuous $h$.  
The role of $F$ in \Cref{L:Local2} is played here by $F=F_h=\int h \d \tilde X^{(l-1)}_{R/2}(w_1)\ldots   \d \tilde X^{(1)}_{R/2}(w_{l-1})  $.  
By the a priori bound \Cref{lem: moment for use} we estimate $(\E[F^4])^{1/4}$ by $C \norm h \norm_{\infty}T^{l-1}$. The role of the set $A$ is played by 
$$
A =  \left\{(s,x) :  s\in [0,T],  \str x \str \leq R   \right \}  \quad \cup\quad  \left\{(s,x) :  s\in [0,T],  \str x-Y_{s} \str \leq R   \right \}  
$$
Note that $A$ is random here, so we actually  need a straightforward generalization of  \Cref{L:Local2} which is omitted.
The probability $\mathbf{P}^{*}_{\nu_k}(E_A)$ is the probability that at least one of the discrepancies started at $\pm(k+1)$ comes closer than $R/2$ to the tagged particle in $[0,T]$ or that a discrepancy enters the region $[-{R/2},{R/2}]$ in time $[0,T]$.  For the discrepancy started at $k+1$, we simply use \Cref{P:Exptag} to argue that the tagged particle can not catch up. For the discrepancy started at $-k-1$, we use the maximal speed of discrepancy and tagged particle, see \Cref{lem: max v dis,,lem: max v tag}. In particular, if $k > R$, and recalling $R>T^{\kappa}, \kappa >1$, we get
\[
\mathbf{P}^*_{\nu_k}[E_A] \leq C \e^{-c (k/T)^{1/2}}.
\]
Performing the sum $\sum_{k>R}(\mathbf{P}^*_{\nu_k}[E_A])^{1/2}$ we get $C\e^{-cT^{\frac{\ka-1}{2}}}$.
This yields the required bound on the variation. 
\end{proof}

\begin{proof}[Proof of \Cref{lem: loc of h}]

It remains to pass from estimates on $Z-\texttt{P}_R Z$ to estimates on $H-H_R$. This is done by telescoping:
\begin{eqnarray}
\str (H-\texttt{P}_RH)(\sigma)\str  &\leq&  \str (f^{(l)}-\texttt{P}_Rf^{(l)})(\sigma) \str   \str Z(\sigma)\str   + \str \texttt{P}_Rf^{(l)}(\sigma) \str     \str (Z-\texttt{P}_RZ)(\sigma)\str    \\[2mm]
&+& \sum_{\beta}\str (v_{\beta}-\texttt{P}_Rv_{\beta})(\sigma) \str   \str Z(\sigma^\beta)\str    +\sum_{\beta}  \str \texttt{P}_R v_{\beta}(\sigma) \str   \str (Z-\texttt{P}_RZ)(\sigma^\beta)\str. \nonumber
\end{eqnarray}
where the index $\beta= (x,\eta, r) \in \Z \times \{1,-1\} \times \N$ and we have abbreviated
$$
v_{\beta} =  \chi^{-\eta}_{r} \chi^\eta_{(x-r,x]} \tilde g^{(l)}_{x,\eta}, \qquad \sigma^\beta= \sigma^{x-r,x}.
$$
We bound only the third term explicitly  (the rest are simpler or similar to handle), call it $V$. 
\begin{equation}
\label{E:biEx}
V(\sigma) \leq \sum_{\beta}   \str (v_{\beta}-\texttt{P}_Rv_{\beta})(\sigma) \str     \str Z(\sigma^\beta)\str (\caS_{l-1}(T)).
\end{equation}
By the triangle inequality, then Cauchy Schwarz
\begin{align}
\|V(\caS_{l-1}(T))\|_{L^q} &\leq \sum_{\beta}  \| v_{\beta}-\texttt{P}_Rv_{\beta})(\si)   \str Z(\sigma^\beta)\str (\caS_{l-1}(T))\|_{L^q}  \\
& \leq  \sup_{\beta} \norm Z(\sigma^\beta)\str (\caS_{l-1}(T)) \norm _{2q} \,  \sum_\beta  \norm (v_{\beta}-\texttt{P}_Rv_{\beta}) \norm_{L^{2q}}  
\end{align}
The $L^{2q}$-norm of $Z(\sigma^\beta)\str (\caS_{l-1}(T))$ is bounded  independently of $\beta$ by $CT^{l-1}$ by \Cref{lem: moment for use}, and remaining sum over $\beta$  is bounded by  $\e^{-cR}$ using \Cref{ass: norms,ass: local}.  This yields the claim. 
\end{proof}

\subsection{An Application to the random measure $K$}\label{sec: app to k}
We localize the measure $K$ too, though in a slightly different sense than for $H$. 
Let 
$$
P_{(-\infty, R]} K(\sigma)=  \Bp (K \, \big \str \,  \sigma_x, x \leq R)
$$
\begin{lemma} \label{lem: loc k}
\beq
\str   P_{(-\infty, R]} K- K \str  \leq  C \e^{-cR} L^{\ell-l}
\eeq
\end{lemma}
\begin{proof}
First we note that we can change $K$ into $K_{R/2}$ (i.e.\ replacing $X^{(i)}$ by $X_{R/2}^{(i)}$ at the expense of an error of order $C \e^{-cR} L^{\ell-l} $ in total variation. Just as for $H$, this is an application of \Cref{lem: moment for use}, see \Cref{lem: loc of z}.  It remains then to estimate $\str   P_{(-\infty, R]} J- J \str $ with $J=K_{R/2}$.We remark that, if $J(\sigma)=\E_{\sigma}(F)$ with $F \in \caB(\sigma_s(x), x \in \Z, s \geq 0)$, then 
\begin{eqnarray}
\Bp^0[\str   P_{(-\infty, R]} J- J \str^q]    & \leq    \E_{\mu(R, \infty), 0} [\str F(\sigma^1)-F(\sigma^2)\str^q]
\end{eqnarray}
If moreover $F \in  \caB(\sigma_s(x), x-Y(s)\leq R/2 )$, then 
$$
\E_{\mu(R, \infty), 0} [\str F(\sigma^1)-F(\sigma^2)\str^q]^{1/q} \leq    C \norm \str F \str \norm_{L^{2q}} \Pr(E)^{1/(2q)}
$$
with $E$ the event that under $\E_{\mu(R, \infty), 0}$, the leftmost discrepancy remains a distance $R/2$ to the right of $Y(t)$ for all times. $\Pr(E) \leq C \e^{-cR} $ by \Cref{P:Exptag}.
\end{proof}

\subsection{Proof of \Cref{L:Factor1}}
We recall that we are out to bound the variation of the measure
$$
  \textrm{d} s \, \E_{\Bp}[H(\si; \d w_{1,\ldots,l-1}) K(\sigma_{T^{\kappa}}, \d u_{1,\ldots, u_{\ell-l}}) ].
$$
on the set $E'$ defined at \Cref{E:EEE}.

Let us fix $s>0$ and estimate the restricted measure uniformly in this variable.  In this case (and when restricted to the relevant subspace of $E'$), $H$ is a measure on $S_1:=\caS_{l-1}(T)$ and $K$ is a measure on $S_2:= \caS_{\ell-l}[T^{\kappa},L]$.   We replace $H,K$ by $P_R H, P_{(\infty,R]} K$. These substitutions make an error in the total variation of order $C L^{\ell-l+1}(\e^{-cR}+\e^{-c T^\alpha}) $, see \Cref{lem: loc of h} and \Cref{lem: loc k}.  Then we are down to estimating 
$$\E(P_RH(\sigma_0)  P_{(\infty,R]} K(\sigma_{T^\ka})) -   \E(P_R H(\sigma_0) ) \E(P_{(\infty,R]} K(\sigma_{T^\ka})))  $$
  and this is now in the form of Lemma \ref{lem: decay of correlations}. This ends the proof of \Cref{L:Factor1}.

\section{Finite Dimensional Convergence and Tightness}
\label{S:CLT}
Having established \Cref{L:Factor1}, we are ready to derive the various functional CLTs.  By adding a constant drift, it suffices to consider only quasi-local processes $X$ with $\E(X(t))=0$.
For such $X$, we consider  $X_n(\cdot):=1/\sqrt{n} X(\cdot \: n)$ and we prove first that for a fixed sequence of times $t_1< \dotsc< t_{l}$ the vector
\beq \label{eq: vector finit}
1/\sqrt{n}(X_n( t_1 ), \dotsc, X_n( t_l ))
\eeq
converges weakly to the appropriate multivariate Gaussian, see  \Cref{T:FDD}. Then we argue that the sequence of processes $(X_n)_{n \in \N}$ is tight in the Skorohod space $D([0, 1], \R)$, see 
\Cref{prop: tight}. By standard reasoning, these two results complete the proof of our main result \ref{thm: main}. The rest of this section is hence devoted to the proof of these results.  

%\subsection{Convergence of finite-dimensional distributions}
We first compute the $t\to \infty$ limit of the variance of $(1/\sqrt{t})X(t)$. It is given as $D=D_{X}:= D_{1}+D_{2}$ with
\begin{align*}
&D_{1}:=  \sum_{x,\eta} \lambda_{\eta}\E( (g_{x,\eta}(0))^2)= \sum_{x,\eta}  \lambda_{\eta} \Bp(g_{x,\eta}^2 ),\\
& D_{2}:=  2\int_0^\infty \E[ h(s)h(0)]  \d s,   \qquad h=f+ \sum_{x,\eta} \lambda_\eta g_{\eta,x}
\end{align*}
It is a straightforward consequence of \Cref{P:Exp} or \Cref{P:Exptag} that $D< \infty$.

%\begin{align*}
%&D_{1}:=  \sum_{x,\eta}\Bp(\langle g^2_{x,\eta}),\\
%& D_{2}:=  \wdr{2}\int_0^\infty \E[ \langle g^{(i)}(s),g^{(i)}(0)\rangle+ f^{(i)}(s)f^{(i)}(0)]  \d s,
%\end{align*}
To prove convergence of finite-dimensional distributions, we use the method of moments. 
Let us consider increasing sequences $(a_i)_{i \leq k}, (b_i)_{i\leq k}\in \R^k$ such that $0\leq a_i < b_i < a_{i+1}$   and let
 \[
 \Delta_i X_n=X_n(b_i)-X_n(a_i) = \frac{1}{\sqrt{n}}[X(b_in)-X(a_in)],
 \]
%$\mathbf X_n:= (\Delta X^{(i)}_n)_{i=1}^k$ and 
%\begin{align*}
%&D_{1, i}:=  \E[ \langle g^{(i)}(0),g^{(i)}(0)\rangle],\\
%& D_{2, i}:=  \wdr{2}\int_0^\infty \E[ \langle g^{(i)}(s),g^{(i)}(0)\rangle+ f^{(i)}(s)f^{(i)}(0)]  \d s,\\
%& D_i=D_{X_i}:= D_{1, i}+D_{2, i}
%\end{align*}
%where $f^{(i)}(t), g^{(i)}(t)$ are
   Let $\mathbf{\gamma}:= (\gamma_i)_{i=1}^k \in \R^k$ and let $(N_{i})_{i=1}^k$ be independent mean zero Gaussians with  respective variance  $D[{b_i}-a_i]$. 

\begin{theorem}
\label{T:FDD}
For all $\ell\in \N$ and $\epsilon>0$, there is $C(\epsilon, \ell)>0$ such that
\[
\bigg|\E \left[(\gamma \cdot \Delta X_n)^\ell\right]- \E\left[(\gamma \cdot N)^\ell\right]\bigg| \leq {Cn^{\epsilon-1/2}}.
\]
This implies (method of moments) that the vector \eqref{eq: vector finit} converges in distribution to $$\sqrt{D}(B(t_1), \ldots, B(t_l)), $$ where  $B(t)$ is a standard Brownian motion.
\end{theorem}

The tightness is also in essence a consequence of the above theorem.
\begin{proposition} \label{prop: tight}

The sequence $(X_n)_{n \in \N}$ is tight in $D_{\Omega}([0, 1])$, equipped with the Skorohod topology. 
\end{proposition}
\begin{proof}
We first fix some notation. 
For any $1>\delta >0$, we fix a partition $\caJ(\delta)$  of $[0,1]$ by intervals $I_j$ with lengths between $\delta$ and $2\delta$.  For any interval $I$ we write
$$
w_{X}(I)=\sup_{s,t \in I} \str  X(t)-X(s)\str
$$
Tightness of the  sequence $X_n$ is implied by the following two conditions (see e.g.\ \cite{billingsley2013convergence}
\begin{enumerate}
\item  For any $\kappa,\epsilon >0$, there is a $1>\delta>0$ such that 
$$
\limsup_n  {\bf P} (\max_{I \in \caJ(\delta)} w_{X_n}(I)   \geq \epsilon )  \leq \kappa
$$
\item For any $\eta>0$, there is an $M$ such that 
$$
\sup_n{\bf P} (\sup_{0\leq t \leq 1} \str X_n(t)\str   \geq M )  \leq  \eta
$$
\end{enumerate}
Now, we check these conditions, starting with $1)$.

As in Section \ref{S:DofC}, we denote by $\bar X$  the quasi-local process obtained derived from $X$ by replacing $(f, g)$ by $(|f|, |g|)$.  Then clearly  $\bar X(t)$ is increasing and so
\[
\sup_{0\leq s\leq t \leq u}|X(t)-X(s)|\leq \bar X(u) \quad \forall u \in \R_+.
\]
Since all $X_n$ are stationary quasi-local processes, we then find
\beq \label{eq: to bound}
{\bf P} (\max_{I \in \caJ(\delta)} w_{X_n}(I)    \geq \epsilon )  \leq C \delta^{-1}   {\bf P} ( \bar X_n(2\delta) \geq \epsilon )
\eeq
To bound the probability on the right hand side, we use \Cref{T:FDD} for $\bar X$ and $\ell=4$; 
$$
\E(\bar X^4_n(t)) \leq  C(t^{2} + n^{-1/4})
$$
so that, by \eqref{eq: to bound} and the Markov inequality, we get 
$$
{\bf P} (\max_{I \in \caJ(\delta)} w_{X_n}(I)  \geq \epsilon )  \leq  C \delta^{-1}  \frac{(\delta^{2} + n^{-1/4})}{\epsilon^4}
$$
which settles condition 1) above.  To check condition 2), it suffices to again consider the increasing $\bar X_n(t)$ and to establish $\E((\bar X_n(1))^2) <C$. The latter follows again by  \Cref{T:FDD}.

\end{proof}

\begin{proof}[Proof of \Cref{T:FDD}]
Much of the work done here is (standard) combinatorics to suitably reduce (by expanding) the moments to expressions we can more easily compute.  Let us fix the time scale $n$.  For simplicity, we first do the case $k=1$. We set $L=b_1-a_1$ and by stationarity we can restrict to the interval $[0,L]$.  We start from
\beq
\label{E:Part1}
\E\left[\left( \gamma \Delta X_n\right)^\ell\right]=n^{-\ell/2} \gamma^{\ell} \int_{[0,L]^\ell} \E\left[\prod_{j=1}^{\ell} \textrm d X (t_j)\right]
\eeq
The measure $\E\left[\prod_{j=1}^{\ell} \textrm d X(t_j)\right]$ is not absolutely continuous due to singular contributions on diagonals  $t_{j}=t_{j'}$.   Formally, this comes about because the powers $(\textrm{d} {X}(t))^{q}$ are not necessarily zero.
We find it computationally convenient to further reduce the problem to a sum of stochastic integrals over the {open simpleces $ \caS_{r}(nL)$, with $r \leq \ell$} by viewing the $(\textrm{d} {X}(t))^{q}$ as quasi-local processes themselves.  

To make this part of the expansion explicit, we introduce more notation:  
Let ${\bf j}=\{j(1), \dotsc, j(r)\}$ denote a partition of $\ell$, i.e.\  $j(l) \in \N$ and $\sum_{l=1}^{r} j(l)=\ell$.
Then 
\beq
\label{E:Part2}
\int_{[0,L]^\ell} \E\left[\prod_{j=1}^{\ell} \textrm d X (t_j)\right]=  \sum_{\bf j}  \frac{\ell !  }{ j(1)! \ldots j(r)!}  \E\left[  Z({{\bf j}})\right]
\eeq
where
%\[
%Z({{\bf j}})= \int_{\caS_{r}(n(b-a))}  \,    \prod_{l=1}^{r} (\textrm d {X}(t_l))^{j(l)}.
%\]
%
%%Let
%%\[
%%Z_i({{\bf j}_{m_i}})= \int_{0 \leq t_l \leq (b_{i}-a_i) n, t_l \neq t_{l'} \forall l, l' \in [r_i]} \prod_{l=1}^{r_i} \textrm d {X^{(i)}}^{j_{m_i}(l)}(t_l+ a_i).
%%\]
%By relabeling, this term can be written in the form
\begin{equation}
\label{E:Part3}
  Z({{\bf j}})= \int_{\caS_{r}(nL)} \textrm{d} W^{(1)}(t_1) \cdots \textrm{d} W^{(r)}(t_r),\qquad   \text{with} \quad \textrm{d} W^{(l)}(t):= (\textrm d {X}(t))^{j(l)} 
\end{equation}
The $W^{(l)}$'s may have increments with nonzero mean.  To give a clean statement below, let $\textrm{d} \caW^{(l)}(t)=  \textrm{d} {W}^{(l)}(t)- \E[ \textrm{d} {W}^{(l)}(t) ] $.  Expanding \cref{E:Part3} gives
\begin{equation}
\label{E:Part4}
\E[ \textrm{d} {W}^{(1)}(t_1) \cdots \textrm{d} {W}^{(r)}(t_r)]= \sum_{A \subset [r]} \prod_{l'\in A^c}\E\left[ \textrm{d} {W}^{(l')}(t_{l'})\right]\E\left[\prod_{l\in A}\textrm{d} \caW^{(l)}(t_l)\right]
\end{equation}
Now we will use input from the previous sections, in particular \Cref{L:Factor1}, to calculate the leading contribution to $ \E\left[\prod_{l\in A} \textrm{d} \caW^{(l)}(t_l)\right]$. 
\begin{lemma}[Iterative Decomposition of Correlations]
\label{L:Iter}
Fix $L\in (0, 1), n \in \N$
Let $(\caW^{(l)})_{l=1}^r$ be quasi-local observables having mean zero increments and respective integrands $(f^{(l)}, g^{(l)})_{l=1}^r$.  Let $r$ be odd, then
\beq
\left\str \, \E\left[\prod_{l=1}^r\textrm{d}\caW^{(l)}(t_l)\right]  \, \right\str_{\caS_r(bn)}  \leq   C n^{r/2-1/2+\epsilon},  \qquad \text{for any $\epsilon>0$},
\eeq
where  $\str \cdot \str_{\caS_r(nL)}$ is the total variation on the simplex $\caS_r(nL)$.  For even $r$, we have
\beq \label{eq: even r}
\left\str \,  \E\left[\prod_{l=1}^r\textrm{d}\caW^{(l)}(t_l)\right] - \prod_{l<r, l \text{odd}}^r \E\left[\textrm{d}\caW^{(l)}(t_l) \textrm{d}\caW^{(l)}(t_{l+1})\right]  \, \right\str_{\caS_r(bn)}   \leq   C n^{r/2-1+\epsilon}. \eeq
Furthermore, the total variation of  $\E\left[\textrm{d}\caW^{(l)}(0) \textrm{d}\caW^{(l+1)}(t)\right]$ on $\{t \geq T\}$ is bounded by $C\e^{-c T^{1/4}}$.
%and the total variation of the second term on the left hand side is bounded by $C(r,\epsilon)(\log^C n)n^{r/2}$.
\end{lemma}

\noindent
This lemma will be proved after the proof of \Cref{T:FDD} is concluded.

We are now ready to determine, from among the terms expanded in \Cref{E:Part1,,E:Part2,,E:Part3,,E:Part4}, the main contributions to the $\ell$'th moment. We keep $\boldsymbol{m}$ fixed and we  compute the contribution from the relevant $\mathbf j$'s and $A$'s.  From \Cref{L:Iter} we deduce that the maximal contribution to \Cref{E:Part4} is of order 
$$
n^{\str A^c\str} n^{\str A \str/2}, 
$$
provided that $i)$ $\str A \str$ is even, $ii)$ for any odd $l \in A$, there  is no $l' \in A^c$ such that $t_l \leq t_{l'} \leq t_{l+1}$, $iii)$ for all $l' \in A$, the increment $\d W^{(l')}$ has nonzero mean. 
Subleading contributions are down by at least a factor $n^{-1/2+\epsilon}$.
%%provided that the paired times times $(t_i, t_{l+1})$ of $A$ do not have any times of $A^c$ between them (ln which case that time would not contribute a factor $n$ by the last statement of \Cref{L:Iter}. 
%All other contributions are down by at least a factor $n^{-1/2+\epsilon}$.  
Looking back at \Cref{E:Part3} and recalling that $\d X$ had zero mean, we see that the leading contributions are of order $n^{\ell/2}$, for $\ell$ even,  and they occur when all $j(l)$ are either $1$ or $2$, and for each time $t_{l}$ for which $j(l)=1$, there is a partner time $t_{l'}$ such that $j(l')=1$ and $\str l-l'\str=1$. The pairs $(l,l')$ are those that constitute the sets $A$ for the dominant contributions in \Cref{L:Iter}.   
Let $q=|\{l: j(l)=2\}|$.  Then the above considerations lead to
\[
\E\left[ Z({{\bf j}})\right]=  \frac{(nL)^{\ell/2} }{(\ell/2)!} D_{1}^{q}(D_{2}/2)^{\ell/2-q} + O(n^{\ell/2-1/2+\eps}).
\]
where we also used that 
$$
\int_{\caS_2(T)} E\left[\textrm{d}X^{}(t_1) \textrm{d}X^{}(t_{2})\right]  = \frac 12 D_{2} T + \caO(C \e^{- cT^{1/4} }), \qquad  \int_{0}^T E\left[ (\d X^{}(t))^2 \right]  = D_{1} T 
$$
%Next, in the cases of interest
%\[
%\pi(m_i; {{\bf j}_{m_i}})= 2^{-q} {m_i \choose q_i}.
%\] 
%Taking into account distinct orderings of times which contribute to $\E\left[\prod_{i} Z_i({{\bf j}_{m_i}})\right]$ gives a further combinatorial factor $\frac{(m-2q_i)!}{(m_i/2-q_i)!} (m_i/2)!$.  The first factor counts the number of ordered pairings of times $(t_i: i\in A)$ and the second factor counts number of ways of choosing $A$ so that the pairs in $A$ are not interlaced with the $(t_i: i\in A^c)$.   We arrive at
%\[
%\E\left[ \prod_{i} Z_i({{\bf j}_{m_i}})\right]=\prod_{i=1}^k 2^{-q_i}m_i!! {m_i/2 \choose q_i}D_{1,i}^{m_i/2-q_i}D_{2,i}^{q_i}(b_i-a_i)^{m_i/2} n^{m_i/2} + O(n^{\ell/2-1+\eps})
%\]
After summing over leading ${\bf j}$ in \eqref{E:Part2},  we arrive at (for $\ell$ even, otherwise we get only the error term)
\[
\E\left[\left( \gamma \Delta X_n\right)^\ell\right]=      (\gamma^2LD)^{\ell/2} \frac{ l! }{(l/2)! 2^{\ell/2}}+O(n^{\ell/2-1/2+\eps})
\]
Recognizing the $\ell$'th moment of a Gaussian on the right hand side concludes the proof for the case $k=1$. For general $k$, we proceed similarly, but with obvious restrictions on the range of time-arguments in the $\d X(t_j)$. The only change that deserves a comment is the case where in \eqref{eq: even r}, there appear pairs $t_l,t_{l+1}$ such that one of them belongs to $[a_i,b_i)$ and the other to 
$[a_{i'},b_{i'})$ with $i\neq i'$. Let us pretend that $a_{i'}=b_i$ (other possibilities are easier to handle). Contributions of such pairs are subleading by the decay of correlation function (last claim of  \Cref{L:Iter}).
 
\end{proof}
\begin{proof}[Proof of \Cref{L:Iter}]
Define first the sequence of numbers $v_i, i=1,2,\ldots$ recursively by 
$$
v_1=\log^4 n,\qquad v_{i+1}= \big(\sum_{j=1}^i v_i\big)^2
$$
The main idea is to decompose the simplex $\caS_r(nL)$ in clusters by grouping consecutive times.  We fix an increasing sequence $(t_1,\ldots,t_r) \in  \caS_r(nL) $ and we define a grouping of the times $t_i$ in clusters (in fact, this is simply a grouping of the indices $1,\ldots, r$).   We let
$\caT_1 :=\{t_1,t_2,\ldots, t_{z_1}\}$ where $z_1>1$ is the first index for which 
$$
(t_{z_1+1}-t_{z_1}) > v_{z_1},\qquad \text{or}   z_1=r
$$
Once $\caT_1$ defined (and $z_1 \neq r$), we define $\caT_2$ by deleting the times $t_{\caT_1}$ from the sequence $(t_1,\ldots, t_r)$, renumbering the remaining ones, and repeating the above step. More concretely, we set
  $\caT_2 :=\{t_{z_1+1},\ldots, t_{z_2}\}$ where  $z_2>z_1$ is the first index for which 
$$
(t_{z_2+1}-t_{z_2}) > v_{z_2-z_1},\qquad \text{or}   z_2=r
$$
This is repeated until we get a cluster decomposition $$\boldsymbol{\caT}= (\caT_1,\ldots, \caT_d)$$ of $(t_1,t_2,\ldots, t_r)$ (some clusters can be singletons). 
  The (sole) important properties of this cluster decomposition are
\begin{itemize}
\item[$a)$]  All times in a given cluster are close to each other: $\max \caT_i-\min \caT_i \leq C\log^{C} n$, where $C=C(r)$
\item[$b)$]  The distance from $\caT_{i+1}$ to $\caT_i$ is large compared to the length of $\caT_i$: There is some $T \geq \log^4n$, depending only on the number of times in $\caT_i$,  such that 
$$
\max \caT_i-\min \caT_i  < T,\qquad    (\min \caT_{i+1} -  \max \caT_i) \geq T^2
$$
\end{itemize}
A cluster decomposition $\boldsymbol{\caT}= (\caT_1,\ldots, \caT_d)$ defines naturally a subset of $\caS_r({nL})$ that we call $\caS_{\boldsymbol{\caT}}$.
We now draw two conclusions from previous estimates, that follow from these respective properties
\begin{enumerate}
\item  The total variation of the measure $\E\left[\prod_{l=1}^r\textrm{d}\caW^{(l)}(t_l)\right] $ on $\caS_{\boldsymbol{\caT}}$ is bounded by $Cn^{d} \log^C n $ with $d$ the number of clusters in the cluster decomposition $\boldsymbol{\caT}$.    This is a consequence of the a-priori estimate of \Cref{lem: moment for use}.
\item  The measure factorizes on clusters, up to a small error:
$$
\left\str \, \E\left[\prod_{l=1}^r\textrm{d}\caW^{(l)}(t_l)\right]  -    \prod_{j}  \E\left[\prod_{t_l\in \caT_j} \textrm{d}\caW^{(l)}(t_l)\right] \,   \right\str_{\caS_{\boldsymbol{\caT}}}  \leq  C n^r \e^{-\log^2 n}
$$
This follows inductively from the crucial \Cref{L:Factor1} by using the property $b)$ above. 
Indeed, by direct application of \Cref{L:Factor1} with $\kappa=2$, we get that
$$
\left\str \, \E\left[\prod_{l=1}^r\textrm{d}\caW^{(l)}(t_l)\right]  -     \E\left[\prod_{t_l\in \caT_1} \textrm{d}\caW^{(l)}(t_l)\right]  \times    \E\left[\prod_{t_l\in \cup_{j\geq 2}\caT_j} \textrm{d}\caW^{(l)}(t_l)\right] \right\str_{\caS_{\boldsymbol{\caT}}}   \leq    C n^{r}  \e^{-\log^2 n}.  
$$
and this is repeated until we have split of all clusters. 
\end{enumerate}

Combining the conclusions 1) and 2), we see that, in total variation 
$$
\E\left[\prod_{l=1}^r\textrm{d}\caW^{(l)}(t_l)\right]  =   \sum_{\boldsymbol{\caT}: d(\boldsymbol{\caT})\geq r/2}    \prod_{j}  \E\left[\prod_{t_l\in \caT_j} \textrm{d}\caW^{(l)}(t_l)\right]  + \caO(n^{r/2-1/2+\epsilon}) 
$$
where $d(\boldsymbol{\caT})$ is the number of clusters. 
However, since the $\textrm{d}\caW^{(l)}(t_l)$ have mean zero, all clusters decompositions $\boldsymbol{\caT}$ with singletons vanish on the right hand side. Hence the only leading contributions are those where each cluster consists of a pair, which proves the lemma except for the last statement. That last statement however follows directly from \Cref{L:Factor1}.
\end{proof}

\bibliography{ToommodelCLT}

\begin{thebibliography}{1}

\bibitem{billingsley2013convergence}
Patrick Billingsley.
\newblock {\em Convergence of probability measures}.
\newblock John Wiley \& Sons, 2013.

\bibitem{borodin2007large}
Alexei Borodin, Patrik Ferrari, and Tomohiro Sasamoto.
\newblock Large time asymptotics of growth models on space-like paths i:
  Pushasep.
\newblock Technical report, WIAS, 2007.

\bibitem{corwin2015q}
Ivan Corwin and Leonid Petrov.
\newblock The q-pushasep: A new integrable model for traffic in 1+ 1 dimension.
\newblock {\em Journal of Statistical Physics}, 160(4):1005--1026, 2015.

\bibitem{crawford2015toom}
Nick Crawford, Gady Kozma, and Wojciech De~Roeck.
\newblock The toom interface via coupling.
\newblock {\em arXiv preprint arXiv:1501.04746}, 2015.

\bibitem{RWRW}
M.~Hil{\'a}rio, F.~den Hollander, V.~Sidoravicius, R.~Soares~dos Santos, and
  A.~Teixeira.
\newblock Random walk on random walks.
\newblock {\em Electronic Journal of Probability}, 20(95):35, 2015.

\bibitem{kipnis1986central}
Claude Kipnis.
\newblock Central limit theorems for infinite series of queues and applications
  to simple exclusion.
\newblock {\em The Annals of Probability}, pages 397--408, 1986.

\bibitem{Kipnis}
Claude Kipnis.
\newblock Central limit theorems for infinite series of queues and applications
  to simple exclusion.
\newblock {\em Ann. Probab.}, 14:273--285, 1986.

\bibitem{sethuraman2000diffusive}
Sunder Sethuraman, SRS Varadhan, and Horng-Tzer Yau.
\newblock Diffusive limit of a tagged particle in asymmetric simple exclusion
  processes.
\newblock {\em Communications on Pure and Applied Mathematics},
  53(8):972--1006, 2000.

\bibitem{Varadhan}
S.~R.~S. Varadhan.
\newblock Self diffusion of a tagged particle in equilibrium for asymmetric
  mean zero random walks with simple exclusion.
\newblock {\em Ann. Inst. H. Poincare (Probabilites)}, 31:273--285, 1995.

\end{thebibliography}
\bibliographystyle{plain}

\flushleft
\begin{tabular}{lr}
\begin{tabular}{l}
Nicholas Crawford \\
Dept. of Mathematics\\
The Technion \\
{\tt njc860@gmail.com}
\end{tabular}
&
\begin{tabular}{l}
Wojciech De Roeck \\
Dept. of Physics\\ 
KU Leuven \\
{\tt wmderoeck@gmail.com}
\end{tabular}
\end{tabular}

%\begin{thebibliography}{10}
%
%\bibitem{billingsley}
%Billingsley, ...
%\bibitem
%A. Borodin and P. Ferrari.  \textit{Large time asymptotics of growth models on space-like paths I: PushASEP},  Electron. J. Probab. 13 (2008), 1380-1418.
%
%
%\bibitem{Cohn}
%Donald L. Cohn, 
%Measure Theory: Second Edition (BirkhŠuser Advanced Texts Basler LehrbŸcher),  July 14, 2013
%
%
%\bibitem{Toom1}
%Crawford, De Roeck, Kozma \textit{The Toom interface via coupling}
%
%
%\bibitem{DLSS}
%B. Derrida, J. L. Lebowitz, E. R. Speer, and H. Spohn. \textit{Dynamics of an anchored Toom interface.}
%J. Phys. A, 24(20):4805--4834, 1991.
%
%\bibitem{LiggettBook}
%T. Liggett. \textit{Interacting Particle Systems.}
% Springer Berlin Heidelberg, (1985).
%\end{thebibliography}

\end{document}